\documentclass[smallextended,referee,envcountsect,]{svjour3}
\smartqed
\usepackage{graphicx}
\usepackage{subfig}
\usepackage{adjustbox}
\usepackage{algorithmic}
\usepackage{algorithm}

\usepackage{amsmath}

\usepackage{amssymb}

\usepackage{multicol}
\usepackage{multirow}

\DeclareMathOperator{\rank}{rank}
\DeclareMathOperator{\card}{card}
\DeclareMathOperator{\Diag}{Diag}
\DeclareMathOperator{\sgn}{sgn}
\DeclareMathOperator{\Prox}{Prox}
\DeclareMathOperator{\argmin}{argmin}

\newcounter{relctr} %% <- counter for relations
\everydisplay\expandafter{\the\everydisplay\setcounter{relctr}{0}} %% <- reset every eq
\newcommand\labelrel[2]{%
	\begingroup
	\refstepcounter{relctr}%
	\stackrel{\textnormal{(\alph{relctr})}}{\mathstrut{#1}}%
	\originallabel{#2}%
	\endgroup
}
\AtBeginDocument{\let\originallabel\label} 

\journalname{Noname}

\begin{document}

\title{Harnessing the mathematics of matrix decomposition to solve planted and maximum clique problem}%Large scale planted clique via matrix decomposition}

\author{Salma Omer \and  M Montaz Ali}
\institute{S Omer,  Corresponding author\at
             University of the Witwatersrand \\
              Johannesburg, South Africa\\
              salmaomer@aims.ac.za
           \and
              M M Ali \at
              University of the Witwatersrand \\
              Johannesburg, South Africa\\
              Montaz.Ali@wits.ac.za
}

\date{Received: date / Accepted: date}
%The correct dates will be entered by the editor.
\maketitle
\begin{abstract}
We consider the problem of identifying a maximum clique in a given graph. We have proposed a mathematical model for this problem. The model resembles the matrix decomposition of the adjacency matrix of a given graph. The objective function of the mathematical model includes a weighted $\ell_{1}$-norm of the sparse matrix of the decomposition, which has an advantage over the known $\ell_{1}-$norm in reducing the error. The use of dynamically changing the weights for the $\ell_{1}$-norm has been motivated. We have used proximal operators within the iterates of the ADMM (alternating direction method of multipliers) algorithm to solve the optimization problem. Convergence of the proposed ADMM algorithm has been provided. The theoretical guarantee of the maximum clique in the form of the low-rank matrix has also been established using the golfing scheme to construct approximate dual certificates. We have constructed conditions that guarantee the recovery and uniqueness of the solution, as well as a tight bound on the dual matrix that validates optimality conditions. Numerical results for planted cliques are presented showing clear advantages of our model when compared with two recent mathematical models. Results are also presented for randomly generated graphs with minimal errors. These errors are found using a formula we have proposed based on the size of the clique. Moreover, we have applied our algorithm to real-world graphs for which cliques have been recovered successfully. The validity of these clique sizes comes from the decomposition of input graph into a rank-one matrix (corresponds to the clique) and a sparse matrix.
\end{abstract}
%All acknowledgements should be placed in the back of the paper after Conclusions..
\section{Introduction}
\label{intro}
Consider an undirected graph $G(V,E)$,  where $V$ is the set of vertices and $E$ is the set of edges, $\vert V\vert=N$. A clique is a subset of vertices of $G$, such that every two distinct vertices are adjacent (a complete subgraph). The maximum clique problem (MCP) is the problem of finding the clique with maximum cardinality from an undirected input graph.
The clique number, $\omega(G)$, on a graph $G$ is the number of vertices in a maximum clique. The MCP belongs to the class of NP-hard problems \cite{karp1972reducibility}. Application areas for the MCP include data mining \cite{eblen2012maximum}, community detection \cite{pattabiraman2015fast,arias2014community}, and bioinformatics \cite{strickland2005optimal,malod2010maximum}.

There are several different formulations of the MCP found in the literature such as the continuous quadratic formulation \cite{hungerford2019general}, integer programming formulation \cite{bomze1999maximum}, and quadratic zero-one formulation \cite{pardalos1992branch}. The indefinite quadratic models have also been suggested \cite{pardalos1990global,pardalos1987constrained,al1990constrained}. 

Tomita and Kameda \cite{tomita2007efficient} have presented a branch and bound algorithm for the MCP based on approximate colouring and proper sorting of the vertices.
There are algorithms suggested for the continuous formulation \cite{belachew2017solving}. Algorithms have also been suggested basing on the mathematical model for the maximum independent set problem \cite{howbert2007maximum}. 

In this paper, we consider the planted clique problem, where a single $n$-node clique is planted first and the remaining non-clique edges are inserted independently with a probability. This problem has been studied by a number of authors \cite{alon1998finding,ames2011convex,ames2011nuclear,feige2000finding}.
A polynomial-time algorithm that finds, almost surely, the unique planted clique of size $n\geq a\sqrt{N}$, for sufficiently large constant $a$, in the random graph $G(N,\frac{1}{2})$ has been presented in \cite{alon1998finding}. 
Feige and Krauthgamer \cite{feige2000finding} have proposed an algorithm based on the Lovasz theta function for finding the planted clique of size $n\geq \Omega (\sqrt{N})$ in $G(N,\frac{1}{2})$. Their algorithm also works in the semi-random hidden clique model, in which an adversary can remove edges from the random portion of the graph.

Ames \cite{ames2011convex}, and Ames and Vavasis \cite{ames2011nuclear} have taken a rank minimization approach in modelling the planted clique problem. In addition, Ames \cite{ames2015guaranteed} has established the guaranteed recovery for planted cliques and dense subgraph when the size of the planted clique (the cardinality of the dense subgraph) is known a priori. 

More recently, Bombina and Ames \cite{bombina2020convex} have  studied the problem of identifying the densest subgraph and densest submatrix problems which can be thought of as a generalization of the MCP.

Our contributions in this paper include a mathematical model of the planted and the maximum clique problem and some theoretical results. Our model follows the core concept of the matrix decomposition problem \cite{candes2011robust,chandrasekaran2011rank} but uses different technique approaches. Our model is such that the integer entries of the solution matrix occur as a natural process of optimization. This was shown by implementing the ADMM algorithm for both the planted and `unplanted' maximum clique problems. Theoretical results established include convergence of the ADMM algorithm, and guaranteed recovery of the low-rank and sparse matrices. We have established conditions that guarantee the recovery and the uniqueness of the solution. Moreover, we have derived a tight bound of the dual matrix that certifies the optimality conditions of our proposed model. Our sufficient conditions are closely related to those given in the references \cite{ames2015guaranteed,ames2011nuclear,bombina2020convex} as they also consider the low rank matrix and its dual. However, our sufficient conditions are stronger in that we do not require a number of additional assumptions such as the number of vertices adjacent to clique vertices  (Theorems~2.1 \& 2.2 \cite{ames2015guaranteed}, also assumptions 1, 2, and 3 in \cite{ames2011nuclear}).%Our sufficient conditions are closely related to those in \cite{ames2015guaranteed,ames2011nuclear,bombina2020convex} as all consider matrices of low rank and they all the sufficient conditions are based on some conditions on the dual matrix of the low-rank matrix. However, Our sufficient conditions are stronger than the ones used in \cite{ames2015guaranteed,ames2011nuclear,bombina2020convex} due to the fact that we do not require a number of additional assumptions such as the number of vertices adjacent to clique vertices in  Theorems~2.1 \& 2.2 \cite{ames2015guaranteed} (also assumptions 1, 2, and 3 in \cite{ames2011nuclear}). 
%due to the fact that we do not require assumptions such as the assumptions in theorems 2.1 and 2.2 in \cite{ames2015guaranteed}, and assumptions 1, 2, and 3 in \cite{ames2011nuclear}, and other assumptions in theorems 2.1 and 3.1 used in \cite{bombina2020convex}.
  Our theoretical and computational results are based on the planted clique problem. However, our algorithm equally works for the general and real-world MCP for which we have provided the numerical evidence.  

The rest of the paper is organized as follows. Section \ref{Sec: Section 2} presents the proposed mathematical model. Section \ref{Sec: Section 3} presents a number of proximal operators which are used in the iterates of the ADMM algorithm. Then we present the convergence of the proposed algorithm in Section \ref{Sec: Section 4}. Section \ref{Sec:section 5n} presents the theoretical recovery and the uniqueness of the solution. In Section \ref{Sec: Section 5}, we evaluate the performance of the proposed algorithm for recovering the planted clique, and also recovering the maximum clique from randomly generated and real-world graphs. Finally, we make concluding remarks in Section \ref{Sec: Section 6}. 

\section{The mathematical model} \label{Sec: Section 2}
We have taken the matrix decomposition approach for the planted clique problem. The matrix decomposition problem separates a given matrix $M$ into its low-rank and sparse component by solving the problem   
\begin{align*}
&\min_{L,S\in \mathbb{R}^{N\times N}}~ \rank(L) + \lambda \Vert S\Vert_0\\
&\qquad \textit{s.t.}~~  L + S =M ,
\end{align*}
where $ \Vert S\Vert_{0} = \card(S)$ is the number of non-zero entries in $S$. Both the rank function and $\ell_0$-norm minimization are non-convex. The nuclear norm $\Vert L\Vert_{*}$ is the sum of singular values $\sigma_i(L)$; it is used as the convex relaxation of rank function, and the $\ell_{1}$-norm, $\Vert S\Vert_{1}=\sum_{i=1}^{N} \sum_{i=1}^{N} \vert S_{ij}\vert$, is used as the convex relaxation of $\Vert S\Vert_{0}$.

In the context of the planted clique problem of size $n$, if we include self-loops and assign 1 to the diagonal elements of the adjacency matrix $M$, then $M$ can be split into a rank-one matrix $L$ (corresponding to the maximum clique) and a sparse matrix $S$.
Thus, the formulation for MCP is 
\begin{align}\label{e11}
&\min_{L,S\in \mathbb{R}^{N\times N}}~ \Vert L\Vert_* + \lambda \Vert S\Vert_1\\
&\qquad \textit{s.t.}~~ M-L-S=0, \label{eq11}\\
&\qquad S_{ij} \in [0,1]\label{eq12},                                    
\end{align}
where $M\in \mathbb{R}^{{N\times N}}$ is the adjacency matrix of graph $G$; $L$ and $S$ are the optimization variables. $\lambda>0$ is the regularization parameter. We refer to problem \eqref{e11}-\eqref{eq12} as the `regular' matrix decomposition formulation for the planted clique problem. 

The main difficulty of the above model is that the entries of the optimal $L$ and $S$ may not be integers. Indeed, this was observed for the nuclear norm minimization model by Ames \cite{ames2011convex} who rounded each entry of the optimal $L$ to the nearest integers. The rounding of entries causes noisy recovery as indicated by a number of figures presented in \cite{ames2011convex}. Below, we propose our mathematical model which can overcome the above difficulty. The main concept of our approach is to use the weighted $\Vert S\Vert_{1}$ norm in the objective function \eqref{e11}. We have demonstrated a posteriori that the weighted $\Vert S\Vert_{1}$ norm is central to achieving integer value of the entries of $L$ and $S$. We have taken a systematic approach to generate the weights.
We approximate each term of $\Vert S\Vert_{0}$, 
\[
\Vert S\Vert_{0}=\#
\begin{cases}
1& \textit{if} ~S_{ij}\neq 0,\\
0 &  ~\textit{otherwise},
\end{cases} 
\] 
 with the function
\begin{align}\label{p1}
\psi(S_{ij})=\frac{S_{ij}}{S_{ij}+\epsilon}, S_{ij}\neq -\epsilon,~ \epsilon>0.
\end{align} 
The function $\psi$ is concave in $\left[ 0,\infty \right) $ with $\psi^{\prime }(S_{ij})>0$ for all $S_{ij}\geq 0$,  $\psi^{\prime}(0)<\infty$.  

\begin{figure}[H]
	\centering
	\includegraphics[width=9.5cm]{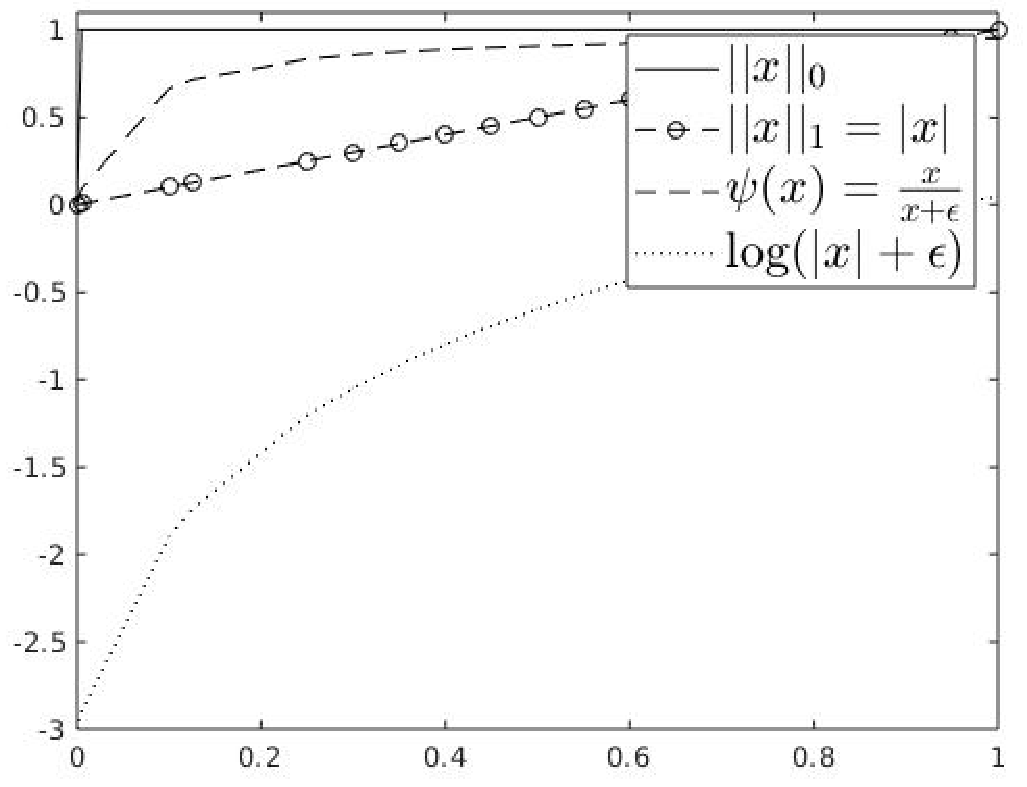}%unInterval-0-1
	\caption{Comparison between $\psi(x)$, the penalty log function $\log(\vert x\vert +\epsilon)$, $\ell_{1}$-norm and $\ell_0$-norm in the case of single variable, $\epsilon=0.05$.} 
	\label{fig:fr1}
\end{figure}

Before presenting the mathematical model we make a graphical comparison of $\psi(x)$ in \eqref{p1} with other, approximations for $\Vert x\Vert_0$ in the case of single variable. The comparison in Figure \ref{fig:fr1} shows that $\psi(x)$, $\epsilon \rightarrow 0$, gives a better approximation of $\Vert x\Vert_0$ than $\vert x\vert $ and $\log(\vert x\vert +\epsilon)$. This motivates our choice of small values of $\epsilon$. 

Hence, our choice of function $\psi(S_{ij})$ further validates the claim made in \cite{candes2008enhancing} that the $\ell_{1}$-norm is not a good approximation for $\ell_0$-norm.
Furthermore, in our mathematical model $S_{ij}=-\epsilon$ can be readily avoided, since $S_{ij} \in [0,1]$.
We write our relaxed objective function as $\Vert L\Vert _*  + \lambda \Phi (S)$,
\begin{align}\label{VARI}
\Phi(S) =\sum_{i=1}^{N} \sum_{j=1}^{N}\psi(S_{ij}).
\end{align}

We now construct a convex surrogate of $\Phi(S)$ at a known feasible $S_{J-1}$, say at $(J-1)$-th iteration of an algorithm. It follows from the concavity of $\psi$ that  

$$\psi(S_{ij})\leq \psi ((S_{J-1})_{ij})  + \psi^{\prime } ((S_{J-1})_{ij}) (S_{ij}-(S_{J-1})_{ij}).$$ 
Hence, we have
\begin{align*}
\Phi(S)& =\sum_{i=1}^{N} \sum_{j=1}^{N} \psi(S_{ij}) \leq \sum_{i=1}^{N}\sum_{j=1}^{N} \left( \psi ((S_{J-1})_{ij})  +  \psi^{\prime }((S_{J-1})_{ij}) (S_{ij}-(S_{J-1})_{ij})\right) \\ & \qquad \qquad\qquad\qquad = \widetilde{\Phi}(S),
\end{align*}
where $\psi^{\prime}((S_{J-1})_{ij})= \frac{\epsilon}{((S_{J-1})_{ij} +\epsilon)^2} $. $\Phi(S)$ satisfies    
\begin{align}\label{eq7}
\Phi(S)\leq \widetilde{\Phi} (S) ~~\textit{and}~~ \Phi (S_{J-1})=\widetilde{\Phi} (S_{J-1}).
\end{align}

\begin{figure}[H]
	\centering
	\includegraphics[width=9.5cm]{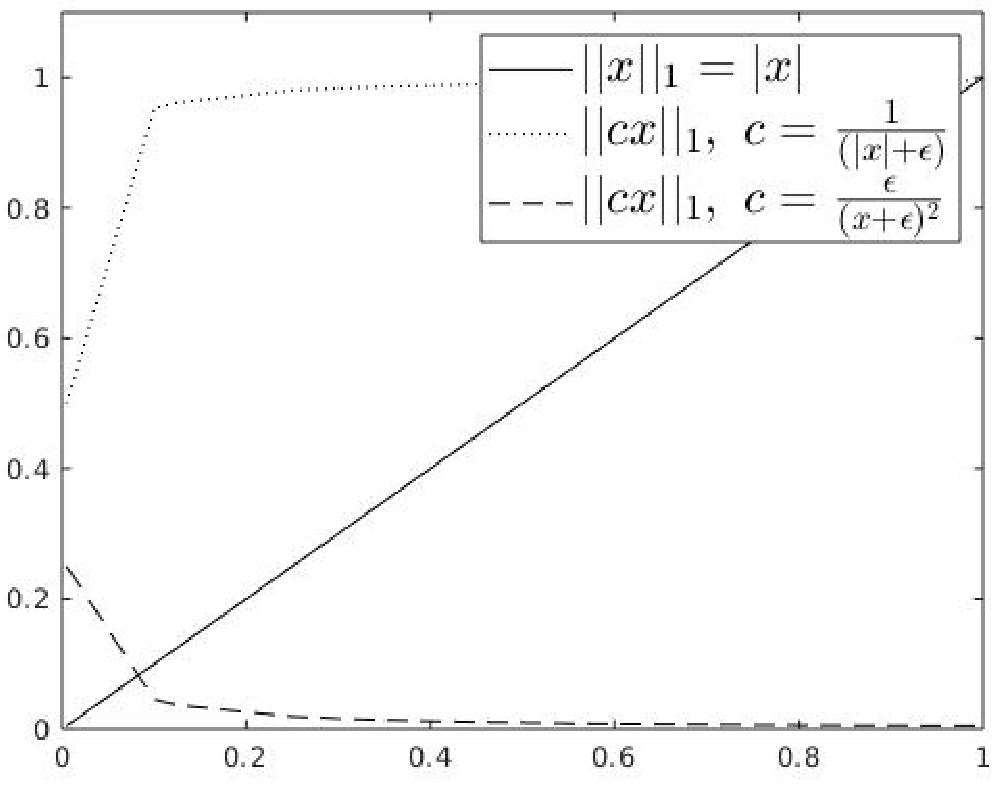}%unInterval-0-1
	\caption{ Comparison between  the convex surrogate function 
		$\Vert cx \Vert_1$, the function $\Vert cx \Vert_1$, $c=\frac{1}{\vert x\vert + \epsilon}$, $\epsilon=0.005$, suggested in in \cite{candes2008enhancing}, and $\ell_1$-norm in the case of single variable.} 
	\label{fig:fr111}
\end{figure}

We now ignore the constant terms in $\widetilde{\Phi} (S)$ and treat the remaining expression as the surrogate for $\Phi(S)$. The concept of the surrogate function has been reported in \cite{han2016two}. Hence, the surrogate becomes $$\sum_{i=1}^{N} \sum_{j=1}^{N}  \psi^{\prime}((S_{J-1})_{ij}) S_{ij}= \sum_{i=1}^{N} \sum_{j=1}^{N} \vert C_{ij} S_{ij}\vert= \Vert C \circ S\Vert_{1},$$ where $C$ is a constant matrix with entries $\frac{\epsilon}{((S_{J-1})_{ij}+\epsilon)^2}$. It is clear that the entries of $C$ are strictly positive.  
The symbol $``{\circ}"$ is known as Hadamard product. Function $ \Phi (S)$ defined in \eqref{VARI} is a concave function. However, the surrogate function $\Vert C \circ S\Vert_{1}$ at $S_{J-1}$ is convex which we refer to as the weighted $\ell_{1}$-norm, where the weights are computed dynamically. We now compare compare our surrogate function, the \lq regular' convex relaxation, and the relaxation suggested in \cite{candes2008enhancing} for the case of single variable in Fig.~\ref{fig:fr111} in the interval $[\epsilon,1]$ for a sufficiently small $\epsilon$. Our surrogate function for single variable case is given by $\Vert cx \Vert_1$, $c= \frac{\epsilon}{({x}+\epsilon)^2}$; $\epsilon=0.005$ is constant.
%We have used a constant value of $c$ for each five data points used to plot $\Vert cx \Vert_1$, e.g. $c_k= \frac{\epsilon}{(\hat{x_k}+\epsilon)^2}$, $\hat{x}_k=x_k$, for data points $x_k, x_{k+1},\cdots,x_{k+5}$, $c_{k+6}= \frac{\epsilon}{(\hat{x_{k+6}}+\epsilon)^2}$, $\hat{x}_{k+6}=x_{k+6}$ for data points $x_{k+6}, x_{k+7},\cdots,x_{k+10}$ and so on. 
%This was done to keep the weight fixed for each epoch which we have discussed next. Fig.2~\ref{fig:fr111} shows the surrogate is a low lying convex in $[\epsilon,1]$. 

% In case of a single constant $c$ and a single variable $x$, the convex surrogate function is given by $\Vert c x \Vert_1 $. 
Figure \ref{fig:fr111} shows the surrogate is a low lying flat like convex in $[\epsilon, 1]$ which allows iterate of an algorithm to land over a range and thereby producing sparse solution via proximal operator.
%As shown in Figure \ref{fig:fr111}, our convex surrogate function $\Vert c x \Vert_1 $, in the interval $[\epsilon,1]$, with the weights $c= \frac{\epsilon}{(x+\epsilon)^2}$, sets most of $x \in [\epsilon,1]$ close to zero better than $\Vert c x\Vert_1$ when the weights are given by $c=\frac{1}{\vert x\vert + \epsilon}$ suggested for the concave function $\log(\vert x\vert +\epsilon)$ in \cite{candes2008enhancing}. %Figure.~\ref{fig:fr111} shows the surrogate is a low lying flat like convex in $[\epsilon,1]$ while allows iterate of an algorithm to land over a range and thereby producing sparse solution via proximal operator.
%This motivates our choice of the convex surrogate function $ \Vert cx\Vert_1$, with weights $c= \frac{\epsilon}{(x+\epsilon)^2}$ extracted from the concave function $\Phi$.
% Comparison between our convex surrogate function $||cx||_1$, $c=\frac{\epsilon}{(x+\epsilon)^2}$, the function $||cx||_1 $, $c=\frac{1}{|x|+\epsilon}$ suggested in \cite{belachew2017solving}  and  $\ell_{1}$-norm in the case of single variable, $\epsilon=0.005$.
 We initialize the optimization algorithm with $(L_0,S_0)$, $C_k=\frac{\epsilon}{((S_{J-1})_{ij}+\epsilon)^2}$, $k=0$, $J=1$, then update the matrix $C$ after every $l$ number of iterations of the algorithm. We then take the corresponding iterates 
\begin{align}\label{iters}
\{(L_0,S_0),(L_{l},S_{l}),(L_{2l},S_{2l}),(L_{3l},S_{3l}),...\},
\end{align}
and denote it as the subsequence $\{(L_k,S_k)\}$, where $(L_k,S_k)$ is the last solution at epoch $k$. The relation between $k$ and $J$ is as follows. If$\mod(J,l)=0$, then we increase $k$ by 1, and update the value of $C_k$ via $C_k =\frac{\epsilon}{((S_{J-1})_{ij}+\epsilon)^2}$, where $\mod(a,b)$ is the modulo operation that finds the reminder after dividing $a$ by $b$. Thus our proposed mathematical model for the planted clique problem and the MCP is given by:  
\begin{align}\label{eq9}
&\min_{L,S\in \mathbb{R}^{N\times N}}~ \Vert L\Vert _*  + \lambda \Vert C \circ S\Vert _1\\
&\qquad \textit{s.t.} ~~\eqref{eq11}-\eqref{eq12},\label{eq91}
\end{align}
where the matrix $C$ is updated at each epoch $k$ of the algorithm used to solve it.

\section{The alternating direction method of multipliers} \label{Sec: Section 3}

In this section, we present the ADMM algorithm which we have used to solve \eqref{eq9}-\eqref{eq91}. We begin with the singular value thresholding and proximal operator which have been used in the iterates of ADMM.

Let $X$ be a matrix of size $N_1\times N_2$ and of rank $r$. Assume that the singular value decomposition (SVD) of $X$ is defined by $X= U\Sigma V^T$, where $U\in \mathbb{R}^{N_1\times r}$, $V\in \mathbb{R}^{N_2\times r}$ and $\Sigma\in \mathbb{R}^{r\times r}$.

For $\tau >0$, we define the SVT (singular value thresholding) operator \cite{cai2010singular} as 
\begin{align}\label{L11}
SVT_{\tau}(X)=U {\mathcal{D}_{\tau}(\Sigma)} V^T,
\end{align}
$\mathcal{D}_{\tau}(\Sigma)=\Diag(\max\lbrace (\sigma_i-\tau ),0\rbrace )$, taking note of the fact that $\sigma_i =\vert \lambda_i \vert $ (eigenvalues) for real symmetric $X$.

Let $ST_{\tau}$ be the  soft thresholding or shrinkage operator with parameter $\tau$. Then for $\tau>0$, $ST_{\tau}$ can be defined on each element of $X$ by:
\begin{align}\label{S11}
ST_{\tau}(X)=\sgn(X)\circ \max(|X|-\tau E,0),
\end{align}
where $\sgn$ is the sign function; $E$ is a matrix of all ones.

The unified shrinkage operator \cite{parikh2014proximal} can be defined by 
\begin{align}
\Prox_{\tau f} (X) = \textit{$\argmin$}_{Y \in \mathbb{R}^{N_1\times N_2}}~ \tau f(Y) +\frac{1}{2} \Vert Y-X\Vert _F^2 \nonumber\\=
\begin{cases}
SVT_\tau (X) & if ~f(Y)= \Vert Y \Vert_{*},\\
ST_\tau(X) & if ~f(Y)= \Vert Y \Vert_{1},
\end{cases} \label{LLA}
\end{align}
where $\Vert .\Vert_{F}$ denotes the Frobenius norm.

It is easy to see that
\begin{align}
\argmin_{Y\in \mathbb{R}^{N_1\times N_2}}& ~~\tau \|  C\circ Y\|_1+\frac{1}{2}\|Y-X \|^2_F\nonumber\\
&= \sgn(X)\circ \max (|X|-\tau C,0)
= STT_{\tau}(C,X)\label{Prox1},
\end{align}
$STT_{\tau}(C,X)=ST_\tau(X)$ when $C=E$.

We now summarize the iterates of the ADMM algorithm for the problem \eqref{eq9}-\eqref{eq91}. 

The augmented Lagrangian of problem \eqref{eq9}-\eqref{eq91} is given by 
\begin{align}\label{AGU}
\mathcal{L}_{\rho}(L,S,\mu)=\|L\|_*+\lambda \|C\circ S\|_1 +\langle\mu, M-L-S \rangle+\frac{\rho}{2}\|M- L-S\|^2_F.
\end{align}
The scaled form of the augmented Lagrangian is as follows:
\begin{align}\label{AGU2}
\mathcal{L}_{\rho}(L,S,\mu)=\|L\|_*+\lambda \|C\circ S\|_1 +\frac{\rho}{2}\|M- L-S+\frac{1}{\rho}\mu\|^2_F- \frac{1}{2\rho}\Vert \mu\Vert_{F}^2,
\end{align}
where $C$ is a constant matrix; $\mu$ is the Lagrange multiplier. The implementation of constants \eqref{eq12} in $\mathcal{L}_{\rho}(L,S,\mu)$ is not required as the iterates of ADMM do not produce negative $S_{ij}$ due to the following reasons. It is straightforward to use an initial feasible solution for problem \eqref{eq9}-\eqref{eq91} solved by ADMM; the entries of input matrix $M$ are $\{0,1\}$; the minimization of $\Vert L\Vert_{*}$ or $\rank (L)$ ensures entries of $L$ cannot be too different.
The ADMM iterates are as follows:
\begin{align*}
L_{J}&=\argmin_{L} ~\mathcal{L}_{\rho}(L,S_{J-1},\mu_{J-1})\\
&=\argmin_L~~  \rho \left( \frac{1}{\rho}\|L\|_*+\frac{1}{2}\|M- L-S_{J-1}+\frac{1}{\rho}\mu_{J-1}\|^2_F\right).
\end{align*}
\begin{align*}
S_{J}&=\argmin_{S} ~\mathcal{L}_{\rho}(L_{J},S,\mu_{J-1})\\
&=\argmin_S~~\rho\left(  \frac{\lambda}{\rho} \|C\circ S\|_1 + \frac{1}{2} \|M-L_{J}-S-\frac{1}{\rho}\mu_{J-1}\|^2_F\right).
\end{align*}
The Lagrangian multipliers $\mu_J$ is updated as follows,
\begin{align}\label{m1}
\mu_{J}=\mu_{J-1} +\rho (M-L_{J}-S_{J}).
\end{align}

Since $\rho $ is a constant, we minimize the following problem using the proximal operator in \eqref{L11} for $L_{J}$:  
\begin{align}
L_{J} &=\argmin_L~~   \frac{1}{\rho}\|L\|_*+\frac{1}{2}\|M- L-S_{J-1}+\frac{1}{\rho}\mu_{J-1}\|^2_F\nonumber \\
&=SVT_{\frac{1}{\rho}}(M-S_{J-1}+\frac{1}{\rho}\mu_{J-1})\label{LK1}.
\end{align}

Similarly, the proximal operator \eqref{Prox1} is used in finding $S_{J}$:
\begin{align}
S_{J}&=\argmin_S~~  \frac{\lambda}{\rho} \|C\circ S\|_1 + \frac{1}{2} \|M-L_{J}-S-\frac{1}{\rho}\mu_{J-1}\|^2_F \label{S333}\\
&=STT_{\frac{\lambda}{\rho}}(C,M-L_{J}+\frac{1}{\rho}\mu_{{J-1}})\nonumber
\\&=\sgn\left( \left[ M-L_J +\frac{1}{\rho} \mu_{{J-1}}\right] _{ij} \right) \nonumber\\
& \qquad \qquad \qquad \quad \circ \max \left(\left | \left[ M-L_J +\frac{1}{\rho} \mu_{{J-1}}\right] _{ij}\right |-\frac{\lambda}{\rho}(C_{J-1})_{ij} ,0\right) \label{3333}.
\end{align}

With the above calculations of the variables $L,~S$, and $\mu$, the steps of the ADMM algorithm of problem \eqref{eq9}-\eqref{eq91} are summarized in \textbf{Algorithm} \ref{a4}.
\begin{algorithm}[H]
	\caption{The ADMM algorithm for problem \eqref{eq9}-\eqref{eq91}}
	\label{a4}
	\begin{algorithmic}
		\STATE{\textbf{Input}: Adjacency matrix $M$, regularization parameter $\lambda>0$, $\rho>0$, $\epsilon>0$, and integer $l$}
		\STATE{\textbf{Initialization}: Start with feasible $L_k=L_{J-1}$, $S_k=S_{J-1}$, $C=C_k (=C_{J-1})$, $\mu_{J-1} =0$, $J=1$, and $k=0$}
		\begin{itemize}
			\item Set $ (C_{J-1})_{ij} = \frac{\epsilon}{( (S_{J-1})_{ij}  +\epsilon)^2}  $, for $i,j\in \left\lbrace 1,2,...,N \right\rbrace $, since $(C_{J-1})_{ij} \in \partial \psi ((S_{J-1})_{ij})$, where $(S_{J-1})_{ij}\in \{0,1\}$
		\end{itemize}
		\*\textbf{WHILE}{ Stopping condition not satisfied}
		\STATE{ \begin{itemize}
				\item $L_{J}=SVT_{\frac{1}{\rho}}(M-S_{J-1}+\frac{1}{\rho}\mu_{J-1})$
				\item $S_{J}=STT_{\frac{\lambda}{\rho}}(C,M-L_{J}+\frac{1}{\rho}\mu_{{J-1}})$
				\item Update $\mu_{J}$ via \eqref{m1},  set $J=J+1$ 
				\item If $\mod(J,l) =0$, then update $ (C_{J})_{ij} =  \frac{\epsilon}{( (S_{J})_{ij}  +\epsilon)^2} $, set $L_k=L_J$, $S_k=S_J$, $C_k=C_J$, $\mu_k=\mu_J$, and $k=k+1$.  
		\end{itemize} }
		\*\textbf{ENDWHILE}
		\STATE{\textbf{Output}: The recovered matrix $L_J$ represents the planted clique}
	\end{algorithmic}
\end{algorithm}

The `regular' model \eqref{e11}-\eqref{eq12} can be solved by adapting \textbf{Algorithm} \ref{a4} where the ADMM iterates are as follows:
\begin{equation}
L_J=SVT_{\frac{1}{\rho}} (M-S_{J-1}+\frac{1}{\rho}\mu_{J-1})\label{L1}, ~\textit{and}
\end{equation}
\begin{align}
S_{J}&=ST_{\frac{\lambda}{\rho}} (M-L_{J}+\frac{1}{\rho}\mu_{J-1})\label{S1}\nonumber\\
&=\sgn\left( \left[ M-L_J +\frac{1}{\rho} \mu_{{J-1}}\right] _{ij} \right) \max \left( \left | \left[ M-L_J +\frac{1}{\rho} \mu_{{J-1}}\right]_{ij}\right |-\frac{\lambda}{\rho} ,0\right).
\end{align}
The update of $\mu_{J}$ is the same as in \eqref{m1}.

 We compare $\lambda \Vert C \circ S\Vert_{1}$ in our model \eqref{eq9}-\eqref{eq91} with $\lambda \Vert S\Vert_{1}$ in the `regular' matrix decomposition model \eqref{e11}-\eqref{eq12}, using the iterates of \textbf{Algorithm} \ref{a4}. A comparison of ADMM iterates \eqref{3333} and \eqref{S1} of $S_J$ shows that \eqref{3333} carries additional information from $(J-1)$-th to $J$-th iteration via $(C_{J-1})_{ij}=\frac{\epsilon}{((S_{J-1})_{ij}+\epsilon)^2}$. It follows from $C_{J-1}$ that $(S_{J-1})_{ij} \longrightarrow 0 \Rightarrow (C_{J-1})_{ij} \longrightarrow 1/\epsilon$, $\epsilon<<1$.

Before making further comparisons, we look at the shrinkage parameter $\tau=\frac{\lambda}{\rho}$ in \eqref{3333} and \eqref{S1}. The theoretical value $\lambda=\frac{1}{\sqrt{N}}$ has been suggested in \cite{candes2011robust}. Clearly, $\tau =\frac{\lambda}{\rho}$ is a small fraction provided $\rho>1$ (a value we have implemented). Comparison of the shrinkage operators \eqref{3333} and \eqref{S1} suggests that if $(S_{J-1})_{ij}$ in previous iteration is small or close to zero then $(\frac{\lambda}{\rho}) \left( C_{J-1}\right)_{ij} $ in \eqref{3333} is larger than $\frac{\lambda}{\rho}$ in \eqref{S1}. This implies that
\begin{equation*}
(S_{J})_{ij}=\sgn\left( \left[ M-L_J +\frac{1}{\rho} \mu_{{J-1}}\right] _{ij} \right) \max \left( \left | \left[ M-L_J +\frac{1}{\rho} \mu_{{J-1}}\right]_{ij}\right |-\frac{\lambda}{\rho} ,0\right),
\end{equation*}
has more chance of staying fractional than 
\begin{align*}
(S_{J})_{ij}&=\sgn\left( \left[ M-L_J +\frac{1}{\rho} \mu_{{J-1}}\right] _{ij} \right)\nonumber\\
& \qquad \qquad \qquad \quad \circ \max \left(\left | \left[ M-L_J +\frac{1}{\rho} \mu_{{J-1}}\right] _{ij}\right |-\frac{\lambda}{\rho}(C_{J-1})_{ij} ,0\right),
\end{align*}
as $\frac{\lambda}{\rho}(C_{J-1})_{ij}>\frac{\lambda}{\rho}.$
At later stages of the algorithm when (majority) entries of $S_{J-1}$ approach towards zero at iteration $J-1$ then this information is fed into iteration $k$ via $C_{J-1}$ with $(C_{J-1})_{ij}>1$, $(C_{J-1})_{ij} \longrightarrow \frac{1}{\epsilon}$ when $(S_{J-1})_{ij} \longrightarrow 0 $. This increases the likelihood of
$\left( \left | \left[ M-L_J +\frac{1}{\rho} \mu_{{J-1}}\right] _{ij}\right |-\frac{\lambda}{\rho}(C_{J-1})_{ij} \right)$ in \eqref{3333} being negative, and thus making $(S_{J})_{ij}=0$. The iterate \eqref{S1} of the regular model does not have this feature, and thus $(S_J)_{ij}$ remains a fraction if $\left | \left[ M-L_J +\frac{1}{\rho} \mu_{{J-1}}\right] _{ij}\right | >\frac{\lambda}{\rho}$. On the other hand, $(S_{J-1})_{ij}$ approaching 1 implies $(C_{J-1})_{ij}<1$. However, in this case, the integer value of $(S_J)_{ij}$ is not an immediate event but rather a gradual optimization process. 
\section{Convergence analysis} \label{Sec: Section 4}
We now present the convergence of the proposed algorithm to the optimal solution. For \textbf{Algorithm} \ref{a4}, we have the following theorem.
\begin{theorem}\label{thm:th1}
	Any accumulation point $(L_*,S_*)$ of the sequence $\left\lbrace (L_{J},S_{J})\right\rbrace $ generated by \textbf{Algorithm} \ref{a4} is an optimal solution of \eqref{eq9}-\eqref{eq91}, with index $J$ corresponding to large epoch $k$.
\end{theorem}
{\it Proof}
\textbf{Algorithm} \ref{a4} computes $(L_{J},S_{J})$ by alternate minimization with respect to one variable while keeping the other one fixed. 
The problem being convex, for sufficiently large $J$ and $\mu_{J-1}$ close to $\mu_*$, and using \eqref{AGU2} we get, 
\begin{align*}
\mathcal{L}_{\rho}(L_{J},S_{J},\mu_{J-1})&= \min_{L,S\in \mathbb{R}^{N\times N }}\mathcal{L}_{\rho} (L,S,\mu_{J-1})\\
&=  \min_{L,S\in \mathbb{R}^{N\times N }}  \Vert L\Vert _*+\lambda \Vert C_J\circ S \Vert _1 + \frac{\rho}{2}\|M- L-S+\frac{1}{\rho}\mu_{J-1}\|^2_F\\&\qquad \qquad \qquad \qquad \qquad \qquad - \frac{1}{2\rho}\Vert \mu_{J-1}\Vert_{F}^2 \\
&\leq \min_{L,S\in \mathbb{R}^{N\times N }}  \Vert L\Vert _*+\lambda \Vert C_J\circ S \Vert _1 + \frac{\rho}{2}\|M- L-S+\frac{1}{\rho}\mu_{J-1}\|^2_F\\
&\leq \min_{L,S\in \mathbb{R}^{N\times N },~L+S=M}  \Vert L\Vert _*+\lambda \Vert C_J\circ S \Vert _1\\&\qquad \qquad \qquad \qquad  \qquad \qquad+ \frac{\rho}{2} \|M-L-S+\frac{1}{\rho}\mu_{J-1} \Vert_{F}^2 \\
&= \Vert L_*\Vert _*+\lambda \Vert C_{*}\circ S_* \Vert _1 +\frac{1}{2\rho} \Vert \mu_{J-1}\Vert_{F}^2,
\end{align*}

since $\left( L_*,S_* \right) $ is the optimal solution of \eqref{eq9}-\eqref{eq91}. It follows that 
\begin{align}\label{las1}
\mathcal{L}_{\rho}(L_{J},S_{J},\mu_{J-1})\leq \Vert L_*\Vert _*+\lambda \Vert C_{*}\circ S_* \Vert _1 + \varepsilon_1,~\varepsilon_1>0,
\end{align}
where $\rho$ is large, $\mu_{J-1}$ is bounded \cite{gao2020admm,magnusson2015convergence} and $\varepsilon_1 = \frac{1}{2\rho} \Vert \mu_{J-1}\Vert_{F}^2 $.

Now it follows from the unscaled Lagrangian in \eqref{AGU}, and \eqref{m1} that
\begin{align}
\Vert L_{J}\Vert _*+\lambda \Vert C_J\circ S_{J} \Vert _1&= \mathcal{L}_{\rho}\left(L_{J},S_{J},\mu_{J-1} \right) - \frac{1}{2\rho}\left(\Vert \mu_{J}\Vert _F^2-\Vert \mu_{J-1}\Vert _F^2  \right)\label{H2}  \\
&\leq  \Vert L_* \Vert_{*} + \lambda \Vert C_{*}\circ S_*\Vert_{1}+\varepsilon_1+ \frac{1}{2\rho} \left(\Vert \mu_{J-1}\Vert _F^2-\Vert \mu_{J}\Vert _F^2  \right) \label{mu1}\\
&\leq  \Vert L_* \Vert_{*} + \lambda \Vert C_{*}\circ S_*\Vert_{1}+ {\varepsilon_1}+\frac{1}{2\rho} \Vert \mu_{J-1} - \mu_{J} \Vert_{F}^2 \label{H3}.
\end{align}
Thus we have
\begin{align}\label{De}
\Vert L_{J}\Vert _*+\lambda \Vert C_J\circ S_{J} \Vert _1 \leq\Vert L_* \Vert_{*} + \lambda \Vert C_{*}\circ S_*\Vert_{1}+ \hat{\varepsilon_1},~ \hat{\varepsilon_1}>0,
\end{align}
where we have used $\hat{\varepsilon_1} = \varepsilon_1+ \frac{1}{2\rho} \Vert \mu_{J-1} - \mu_{J} \Vert_{F}^2 $ using the boundedness of $\{\mu_J\}$ and large $\rho$. The inequality in \eqref{mu1} follows from \eqref{las1} and the inequality \eqref{H3} follows from the reverse triangular inequality applied to $\mu_{J}$ and $\mu_{J-1}$ (reverse triangle inequality holds for any matrix norm \cite{thompson1978matrix,ipsen2009numerical}). 

Using the reverse triangular inequality, $\Vert A-B \Vert_{*}\geq \Vert A\Vert_{*}- \Vert B\Vert_{*}$, of the nuclear norm, and by the optimizer $(L_*,S_*) $ we get
\begin{align}
\Vert L_{J}\Vert _*+\lambda \Vert C_J\circ S_{J} \Vert _1&\geq \Vert M-S_{J}\Vert _*-\Vert M-L_{J}-S_{J}\Vert _*+\lambda \Vert C_J\circ S_{J}\Vert _1 \nonumber \\
&\geq  \Vert L_* \Vert_{*} + \lambda \Vert C_{*}\circ S_*\Vert_{1} -\Vert M-L_{J}-S_{J}\Vert _* \nonumber\\
&=\Vert L_* \Vert_{*} + \lambda \Vert C_{*}\circ S_*\Vert_{1}-\frac{1}{\rho}\left(\Vert \mu_{J}-\mu_{J-1}\Vert _*  \right), \textit{ by } \eqref{m1}, \label{mu2}
\end{align}
where we have used $A= M-S_{J}$, $B=  M-L_{J}-S_{J}$, and assumed $M- S_J$ being the deviation from $L_*$ (noting that $M-S_*=L_*$) and hence $\Vert M-S_J\Vert_* \geq \Vert L_*\Vert_{*}$.
This together with the boundedness of $\{\mu_J\}$ imply
\begin{align}\label{k2}
\Vert L_{J}\Vert _*+\lambda \Vert C_J\circ S_{J} \Vert _1&\geq \Vert L_* \Vert_{*} + \lambda \Vert C_{*}\circ S_*\Vert_{1} -\hat{\varepsilon_2},~ \hat{\varepsilon_2}>0,
\end{align}
where $\hat{\varepsilon_2} = \frac{1}{\rho}\Vert \mu_{J}-\mu_{J-1}\Vert _*$.

Thus, from equations $\eqref{De}$ and $\eqref{k2}$ we have,
$$-\hat{\varepsilon}\leq \Big( \Vert L_{J}\Vert _*+\lambda \Vert C_J\circ S_{J} \Vert _1\Big)  -\Big( \Vert L_* \Vert_{*} + \lambda \Vert C_{*}\circ S_*\Vert_{1}\Big)  \leq \hat{\varepsilon},$$
where $ \hat{\varepsilon} =\max \{\hat{\varepsilon}_1,\hat{\varepsilon}_2\}$. This implies 
$$\Biggl\lvert ~\Big( \Vert L_{J}\Vert _*+\lambda \Vert C_J\circ S_{J} \Vert _1\Big)  -\Big( \Vert L_* \Vert_{*} + \lambda \Vert C_{*}\circ S_*\Vert_{1}\Big) ~\Biggr\rvert\leq \hat{\varepsilon}.$$
For sufficiently large $\rho$, $ \hat{\varepsilon}\longrightarrow 0$ and $ \Vert L_{J}\Vert _*+\lambda \Vert C_J\circ S_{J} \Vert _1$ converges to $ \Vert L_{*}\Vert _*+\lambda \Vert C_{*}\circ S_{*} \Vert _1$. %\Halmos
\qed
In addition, we have established another convergence result as stated in the following theorem by demonstrating that any limit point in an iteration sequence generated by \textbf{Algorithm} \ref{a4} is a KKT point.
\begin{theorem}\label{thm:th2}
	The limit point $(L_*,S_*)$  of the sequence $\{(L_J,S_J)\}$ generated by \textbf{Algorithm} \ref{a4} is the KKT point for problem \eqref{eq9}-\eqref{eq91}.
\end{theorem}
{\it Proof}

We begin by showing the boundedness of $\{L_J\}$ and $\{S_J\}$. It follows that
\begin{align}
\mathcal{L}_{\rho} (L_{J+1},S_{J+1},\mu_J) &\leq \mathcal{L}_{\rho} (L_{J+1},S_{J},\mu_J) \nonumber \\&\leq \mathcal{L}_{\rho} (L_{J},S_{J},\mu_{J})=\mathcal{L}_{\rho} (L_{J},S_{J},\mu_{J-1}) \nonumber \\&\qquad \qquad \qquad \qquad \quad +\langle \mu_J - \mu_{{J-1}},M-L_J-S_J \rangle \label{H4}  \\
& = \mathcal{L}_{\rho} (L_{J},S_{J},\mu_{J-1}) + \frac{1}{\rho} \Vert \mu_{J} -\mu_{J-1} \Vert_{F}^2, \nonumber
\end{align}
where we have used the augmented Lagrangian \eqref{AGU}; the equality in \eqref{H4} follows by writing 
$\mathcal{L}_{\rho} (L_J,S_J,\mu_J)=\Vert L_J\Vert_{*} +\lambda \Vert C_J \circ S_J\Vert_{1} +\langle \mu_{J-1},M-L_J-S_J\rangle +\frac{\rho}{2} \Vert M-L_J-S_J\Vert_{F}^2+\langle \mu_{J},M-L_J-S_J\rangle -\langle \mu_{J-1},M-L_J-S_J\rangle $. Then by recalling the boundedness of $\{\mu_J\}$, $\frac{1}{\rho} \Vert \mu_{J} -\mu_{J-1} \Vert_{F}^2<\infty$, we have 
\begin{align} \label{PR2}
\mathcal{L}_{\rho} (L_{J+1},S_{J+1},\mu_J) - \mathcal{L}_{\rho} (L_{J},S_{J},\mu_{J-1} )<\infty, \forall J.
\end{align}
Thus, it follows from \eqref{las1} and \eqref{PR2} that $\mathcal{L}_{\rho} (L_{J},S_{J},\mu_{J-1}) $ is bounded. On the other hand, from \eqref{H2} we have 
\begin{align*}
\Vert L_J \Vert_{*} + \lambda \Vert C_J \circ S_J \Vert_{1} = \mathcal{L}_{\rho} (L_J,S_J, \mu_{J-1})  -\frac{1}{2\rho} (\Vert \mu_{J} \Vert_{F}^2 -\Vert \mu_{J-1}\Vert_{F}^2),
\end{align*}
for which $ \Vert \mu_{J} \Vert_{F}^2 -\Vert \mu_{J-1}\Vert_{F}^2$ and $\mathcal{L}_{\rho} (L_J,S_J, \mu_{J-1}) $ are bounded. Hence $\Vert L_J \Vert_{*}$ and $ \lambda \Vert C_J \circ S_J \Vert_{1}$ are also bounded. Therefore, both $\{L_J\}$ and $\{S_J\}$ are bounded.

We have shown in Theorem \ref{thm:th1} that $\{(L_J,S_J)\}$ converges to $(L_*,S_*)$ for $J\longrightarrow \infty$. The KKT conditions of problem \eqref{eq9}-\eqref{eq91} are 
\begin{equation}\label{KKTC}
\left\{
\begin{array}{l}
0\in \partial \Vert L_*\Vert_{*}-\mu_*,\\
0\in \lambda \partial \Vert C_{*}\circ S_*\Vert_{1}-\mu_*,\\
M-L_*-S_*=0.

\end{array}
\right.\\
\end{equation}

By the boundedness of $\{\mu_J\}$, $\{L_J\}$ and $\{S_J\}$ we have 
\begin{align}\label{NM}
\lim\limits_{J\longrightarrow \infty} (M-L_J-S_J) =\frac{1}{\rho} \lim\limits_{J\longrightarrow \infty}(\mu_J - \mu_{J-1}) =0,
\end{align}
thus, we have $M-L_*-S_*=0$.

The optimizers $L_J$ and $S_J$ of the sub-problems using \eqref{AGU2} at the $J$-th iteration of ADMM imply 
$$0\in \partial_L \mathcal{L}_{\rho} (L_J,S_J, \mu_{{J-1}}),~~0\in \partial_S \mathcal{L}_{\rho} (L_J,S_J, \mu_{{J-1}}),$$ 
that is 
$$0\in \partial \Vert L_{J}\Vert_{*}-\mu_{J-1}-\rho(M-L_{J}-S_{J}) ,~~ 0\in \lambda \partial  \Vert C_J\circ S_{J}\Vert_{1}-\mu_{J-1}-\rho(M-L_{J}-S_{J}).$$
Hence, there exist $Y_J\in \partial \Vert L_{J}\Vert_{*}$, and $Z_J\in \lambda \partial \Vert C_J\circ S_{J}\Vert_{1}$ such that $$Y_J-\mu_{J-1}-\rho(M-L_{J}-S_{J})=0, ~~Z_J-\mu_{J-1}-\rho(M-L_{J}-S_{J})=0, J\longrightarrow \infty.$$
It follows from \eqref{NM} that $Y_*-\mu_{*}=0,$ and $Z_*-\mu_{*}=0$, $J\longrightarrow \infty$. This implies that 
\begin{equation*}
\left\{
\begin{array}{l}
\mu_*\in \partial \Vert L_*\Vert_{*},~~\mu_*\in \lambda \partial  \Vert C_{*}\circ S_*\Vert_{1},\\
0\in \partial \Vert L_*\Vert_{*}-\mu_{*},~~0\in \lambda \partial  \Vert C_{*}\circ S_*\Vert_{1}-\mu_*.

\end{array}
\right.\\
\end{equation*}
Now we can see that $(L_*,S_*,\mu_*)$ satisfies the KKT conditions \eqref{KKTC}.% \Halmos
\qed

\section{ Theoretical guarantee for exact recovery}\label{Sec:section 5n}

We begin with some preliminaries. Let the rank of symmetric $L_*$ be $r$. Hence, $L_*$ is orthogonally diagonalizable. Then $L_*= U \Sigma U^T$, $U=[u_1,u_2,...,u_r]$, where $u_i$ it the $i$-th singular vector of $L_*$. $\Sigma=\Diag(\sigma_1, \sigma_2,...,\sigma_r)$ where $\sigma_{i}$ is the $i$-th singular value of $L_*$. We assume that the number of non-zero entries in $S_*$ is $m$, i.e., $\vert supp(S_*)\vert = m$. It is easy to see that the support sets of $S_*$ and $C_*\circ S_*$ are equal and hence $\vert supp(S_*)\vert =\vert supp(C_*\circ S_*)\vert $.

Denote by $\mathcal{R}$ the linear space of matrices 
\begin{align*}
\mathcal{R} := \left\lbrace UX^T + YU^T ~\big \vert ~ X,~Y \in \mathbb{R}^{N\times r} \right\rbrace.
\end{align*}

The orthogonal projection $\mathcal{P}_\mathcal{R}$ onto $\mathcal{R}$, is given by:
\begin{align*}
\mathcal{P}_\mathcal{R}(X) = UU^T X + X UU^T - UU^TXUU^T,
\end{align*}
and $\mathcal{P}_\mathcal{R}^\perp (X) = \left( \mathcal{I} -UU^T \right) X \left(\mathcal{I} -UU^T \right) $ is the orthogonal complement projection onto $\mathcal{R}$, where $\mathcal{I}$ is the identity operator. For any matrix $X$, $\Vert \mathcal{P}_\mathcal{R}^\perp X\Vert \leq \Vert X\Vert $ holds, where $\Vert.\Vert$ denotes the spectral norm.

Let us define the linear space of sparse matrices by
\begin{align}\label{eq:4}
\Omega := \left\lbrace S\in \mathbb{R}^{N\times N} ~\big \vert ~ \vert supp(S)\vert =m\right\rbrace. 
\end{align}

Define $\mathcal{P}_{\Omega}$ to be the orthogonal projection onto $\Omega$, that is,
\[\mathcal{P}_\Omega (X) = 
\begin{cases}
X, &\textit{ if }X\in \Omega,\\
0,& \textit{otherwise},
\end{cases}
\]
then $\mathcal{P}_\Omega^\perp$ defined by $\mathcal{P}_\Omega^\perp (X) = X - \mathcal{P}_\Omega (X)$ represents the orthogonal complement projection onto $\Omega$.

The sub-gradient of the $\ell_1$-norm at $C_*\circ S_*$ is of the form $\sgn(C_*\circ S_*)+F$, $\mathcal{P}_{\Omega}(F)=0$, $\Vert F\Vert_{\infty}\leq 1$, where $\Vert F\Vert_{\infty}$ denotes the largest element of $F$ in magnitude. Also the sub-gradient of the nuclear norm at $L_*$ is of the form $UU^T + W$, $\mathcal{P}_{\mathcal{R}}W=0$, $\Vert W\Vert \leq 1$. We will be writing $C\circ S_*$, replacing the constant matrix $C_*$ with $C$. 

We now discuss the regularization parameter $\lambda$ and the rank-sparsity incoherence in the context of our problem. The value of $\lambda$ and the satisfaction of incoherence conditions play central role in the recovery of $(L_*,S_*)$.

The value $\lambda=\frac{1}{\sqrt{N}}$ suggested in \cite{candes2011robust} follows the inverse square root law. We would like to make $\lambda$ dependent on the prior information of the problem at hand. In particular, we use the size $n$ of the clique and define $\lambda $ to be $\lambda= \frac{\frac{N}{m}}{\sqrt{N}}<1$. Then $\lambda$ follows inverse square root law  provided that $\frac{N}{m}<1$. The value of $m$ must obey $ m<(N^2-n^2)$ due to the fact that $n$ is the clique size in $L_*$ and $S_*=M-L_*$. We make a reasonable choice for the size of the sparsity by taking $m= p(N^2-n^2)$; $p = \frac{1}{2}$ is generally used in the planted clique problem \cite{ames2011convex}. We restrict our planted clique size such that $ c_1N\leq n\leq c_2 N$, $n=c{N}$, where $c\in [c_1,c_2]\approx [0.1, 0.9]$. It follows that 
\begin{align*}
\frac{m}{N} = \frac{p\left( N^2-n^2\right) }{N}= \frac{p\left( N^2-c^2N^2\right) }{N} = p\left( {1-c^2}\right)  N > 1, \textit{ for large }N.
\end{align*}
The above inequality holds even for $N=15$, $c=0.9$, and $p=0.1$. Hence ${\frac{N}{m}}<1$ holds. Our choice of $\lambda$ is therefore given by 
\begin{align}\label{Lambda}
\lambda= \frac{\frac{N}{m}}{\sqrt{N}}=\frac{\alpha}{\sqrt{N}}, \textit{ where } \alpha=\frac{N}{m}<1.
\end{align}
We have estimated a range, $[0.0021,0.0914]$, of values of $\alpha$ numerically by plotting $\alpha$ against $N$ for a number of clique sizes $n$ in $m=p(N^2-n^2)$, see Section \ref{Sec: Section 5}.

We now present the incoherence conditions. 
By construction, $M$, $L_*$ and $S_*$ are all symmetric matrices. With this prior information we now present the conditions on $L_*$ and $S_*$ for their guaranteed recovery. We begin with the rank-one matrix $L_*$. It is easy to see that $\sigma_1 =\Vert L_*\Vert_{*}=\sqrt{\sum_{i=1}^{N} \sigma_{i}^2}=\Vert L_*\Vert_{F}=\sqrt{\sum_{i=1}^{N} \sum_{j=1}^{N} (L_*)_{ij}^2} =n$, since $\sigma_{i}=0$, $i=2,3,...,N$. 
It follows that the elements of $UU^T$ are from $\{0,\frac{1}{n}\}$ and the elements of $U$ are from $\{0,\frac{1}{\sqrt{n}}\}$, since $L_*=nUU^T$. The fact that column/row spaces of $L_*$ are not closely aligned with the canonical basis vectors is guaranteed with the following condition proposed in \cite{candes2011robust} 

\[\max_i\Vert U^T e_i\Vert_2^2 \leq \frac{\mu_0 r}{N},~~\forall i, ~~ 1\leq \mu_0 \leq \frac{N}{r}. \]
It is easy to see that for our problem
\begin{align}\label{eq:1}
\frac{1}{n}=\max_i \Vert U^Te_i\Vert_2^2 \leq \frac{\mu_0 r}{N},~~\forall i,
\end{align}
holds for $ 1\leq \mu_0\leq  \frac{N}{r}$.

The joint incoherent condition presented in \cite{candes2011robust} is defined by 
\begin{align}\label{eq:33}
\Vert UV^T\Vert_{\infty}< \sqrt{\frac{\mu_1 r}{N^2}} , ~\mu_1=\mu_0^2 r, ~ 1\leq \mu_0\leq \frac{N}{r}.
\end{align} 
It is also easy to see that the above condition also holds for $1\leq \mu_0 \leq \frac{N}{r}$ for our problem since
\begin{align}\label{eq:3}
\frac{1}{n}=\Vert UU^T\Vert_{\infty}\leq  {\frac{ \mu_0 r}{N}}.
\end{align} 

For the guaranteed recovery of $(L_*,S_*)$ the condition on $S_*$ is that its sparsity pattern is not too structured. This can be achieved by considering Bernoulli model with probability $ \rho$. However $S_*$ must be symmetric and its construction is such that $\vert supp(S_*)\vert<N^2-n^2$; the probability $\rho$ can be adjusted for this support. These properties are needed to ensure feasibility of constraint \eqref{eq11}. Hence, we work with the empirical probability $p$ and the construction of which is as follows. We divide the set of entry locations of $S_*$ into three sets, $A_1=\{(i,j)~\vert ~ i<j\}$, $A_2=\{(i,j)~\vert ~ i>j\}$, and $A_3=\{(i,j)~\vert ~ i=j\}$, $\forall i,j \in \{1,2,...,N\}$. We then apply Bernoulli model in the set $A_1 \cup A_3$. The entry values corresponding to $A_1$ are then copied to set $A_2$. We then calculate the empirical probability $p$ using entries in $S_*=A_1\cup A_2\cup A_3$, treating the entry values as the results of random experiments (Since non-zero entries of $S_*$ are formed using Bernoulli probability model and $M-S_*=L_*$, the planted clique location can be considered random).    

To see the structure of sparsity pattern of $S_*$ we calculate the variance of the elements of each row or column. Each component of a row associates a random variable which assumes 1 with probability $p$ and 0 with probability $1-p$. This implies that the mean and variance of the random variables are $p$ and $p(1-p)$. The expected cardinality of a row or column is $Np$, and similarly $Np(1-p)$ for the variance. Hence we can see that no pattern is guaranteed since the expected value $Np(1-p)$ is the same for every column. That is 
\begin{align}\label{SV}
\big \vert  Var(S_{i*}) - Var(S_{j*}) \big \vert <\varrho,
\end{align} 
for any $\varrho >0$, where $Var(S_{i*})$ is the variance of the entries of the $i$-th column.

Given the above incoherence conditions on $L_*$ and $S_*$, the recovery is guaranteed by the convex optimization. We have the following theorem.

\begin{theorem}\label{thmm:th1}
	Suppose $L_*$ is an $N\times N$ matrix of rank $r$ which obeys incoherence conditions \eqref{eq:1} and \eqref{eq:3}. Moreover, entries for all the rows or columns of $S_*$ satisfy \eqref{SV}. Then there is a numerical constant $c$ such that with probability at least $1-cN^{-10}$, the output $(L_*,S_*)$ of the optimization problem 
	\begin{align}\label{eq:m}
	\min_{L,S}~& \Vert L\Vert_{*} + \lambda \Vert C \circ S\Vert_{1}\\
	\qquad &\textit{ s.t. } M=L+S, \nonumber\\
	&\qquad S_{ij} \in [0,1],\label{eq:m11}
	\end{align} 
	$\lambda = \frac{\alpha}{\sqrt N}$, $0.0021 <\alpha< 0.0914$, is exact, provided that $$\rank(L_*)\leq \tau  \sqrt{\frac{N}{\rho_0 \log^2N}},~ \tau \geq 1,~~~ m< N^2-n^2,$$ 
	$\rho_0 =\frac{m}{N^2}$, $m<N^2-n^2$, is a numerical constant. 
\end{theorem}

The above theorem ensures the decomposition of adjacency matrix $M$ into a  rank-one
matrix, representing the clique, and a sparse matrix.

Based on approximate dual certificates, we establish Lemmas \ref{l2}-\ref{l4}, and the proof of Theorem \ref{thmm:th1} follows from the lemmas.

We establish some conditions for the pair $(L_*,S_*)$ to be the unique optimal solution to our proposed model. These conditions, expressed in terms of the dual matrix $W$, are given in Lemma \ref{l2} which is similar to Lemma 2.4 in \cite{candes2011robust}. However, we have tightened the conditions by using different bounds for our proof. This was possible due to the fact that the conditions $U^T W=0$, $UW=0$ and $\Vert W\Vert<1$ must hold. $U\in \mathbb{R}^{N\times 1}$ has exactly $n$ non-zero entries since $L_* = nUU^T$; each non-zero element equals to $\frac{1}{\sqrt{n}}$. Satisfaction of $U^TW=0$ and $WU=0$ imply that most elements of $W$ must be zero, and the non-zero elements $W_{ij}$ of $W$ must be very small in magnitude so that $\frac{1}{\sqrt{n}} \times W_{ij}\approx 0$, $i,j\in \{1,2,...,N\}$. This results in $\Vert W\Vert <<1$.  Hence, for a moderate approximation of non-zero $W$ we suggest it satisfies $\Vert W \Vert \leq  \frac{\alpha}{2}$, where $\alpha $ is defined in \eqref{Lambda}.  
\begin{lemma}\label{l2}
	Assume that the subspaces $\mathcal{R}$ and $\Omega$ have a trivial intersection, $\Omega \cap \mathcal{R} =\{0\}$, with $ \Vert \mathcal{P}_{\Omega} \mathcal{P}_{\mathcal{R}}\Vert\leq \frac{1}{2}$ and $\lambda <\alpha$, and assume that there exist $W\in \mathcal{R}^\perp$, $F\in \Omega^{\perp}$, and $B\in \Omega$, for which the following conditions hold:
	\begin{enumerate}
		\item  \label{itm:first} $ UU^T+W=\lambda (\sgn(C\circ S_*)+F+\mathcal{P}_\Omega B)$,
		\item  \label{itm:snd} $\mathcal{P}_\mathcal{R}(X)=UU^T$ and $\mathcal{P}_\Omega(Y)=\sgn(C\circ S_*)$, where $X\in \partial \Vert L_*\Vert_{*}$ and $Y\in \partial \Vert C\circ S_* \Vert_{1}$, respectively,
		\item \label{itm:thd} $\Vert W\Vert \leq \frac{\alpha}{2}$, $\Vert F\Vert_{\infty}<\frac{1}{2}$ and $\Vert \mathcal{P}_\Omega B\Vert_F \leq \frac{1}{4}$, for $0.0021 <\alpha< 0.0914$.
	\end{enumerate} 
	($F$ and $B$ will be introduced in Section \ref{Sec: Section 4}). Then $(L_*,S_*)$ is a unique solution to problem \eqref{eq:m}.
\end{lemma}

{\it Proof}
Consider any feasible solution $(\hat{L},\hat{S})=(L_*+P,S_*-P)$ to \eqref{eq:m} such that $(\hat{L} -L_*,\hat{S}-S_*)\neq(0,0)$. It is clear that this feasible solution is a perturbation of the optimal solution $(L_*,S_*)$, and it satisfies the feasibility constraint in \eqref{eq:m}. We show that $G(\hat{L},\hat{S})> G(L_*,S_*)$ for non-zero $P$, where $G(.,.)$ is the objective function in \eqref{eq:m}.

Let $X\in \partial \Vert L_* \Vert_{*}$ and $Y\in \partial \Vert C\circ S_* \Vert_{1}$, then by the definition of the sub-gradient we have 
\begin{align}\label{myeq}
G(\hat{L},\hat{S})&\geq  G(L_*,S_*) + \langle X,\hat{L}-L_*\rangle+\lambda \langle Y,\hat{S}-S_* \rangle \nonumber\\
&=  G(L_*,S_*) + \langle X,\hat{L}-L_*\rangle+\lambda \langle Y,\hat{S}-S_* \rangle \nonumber\\&\qquad \qquad \qquad  \qquad + \langle UU^T+W,\hat{L}-L_* \rangle-\langle UU^T+W,\hat{L}-L_* \rangle \nonumber\\
&\qquad \qquad \qquad \qquad \qquad +\lambda \langle \sgn(C\circ S_*)+F+\mathcal{P}_\Omega B,\hat{S}-S_* \rangle \nonumber\\&\qquad \qquad \qquad \qquad \qquad \qquad  -\lambda \langle \sgn(C\circ S_*)+F+\mathcal{P}_\Omega B,\hat{S}-S_* \rangle \nonumber\\
&\labelrel={myeq:11}G(L_*,S_*) + \langle X-UU^T-W,\hat{L}-L_*\rangle \\&\qquad \qquad \qquad +\lambda \langle Y-\sgn(C\circ S_*)-F-\mathcal{P}_\Omega B,\hat{S}-S_* \rangle \nonumber\\
&\qquad \qquad \qquad \qquad \qquad \qquad  +\langle UU^T+W,\hat{L}-L_*+\hat{S}-S_* \rangle \nonumber\\
&\labelrel={myeq:1} G(L_*,S_*) + \langle X-\mathcal{P}_\mathcal{R}(X)-W,\hat{L}-L_*\rangle \nonumber\\&\qquad \qquad \qquad \qquad +\lambda \langle Y-\mathcal{P}_\Omega (Y)-F-\mathcal{P}_\Omega B,\hat{S}-S_* \rangle \nonumber \\
&\labelrel={myeq:2} G(L_*,S_*) + \langle \mathcal{P}_\mathcal{R}^\perp(X)-W,\hat{L}-L_*\rangle\nonumber\\&\qquad \qquad \qquad \qquad \qquad \qquad \qquad   +\lambda \langle \mathcal{P}_\Omega^\perp (Y)-F-\mathcal{P}_\Omega B,\hat{S}-S_* \rangle \nonumber\\
&= G(L_*,S_*) + \langle \mathcal{P}_\mathcal{R}^\perp(X),\hat{L}-L_*\rangle - \langle W,\hat{L}-L_*\rangle+\lambda \langle \mathcal{P}_\Omega^\perp (Y),\hat{S}-S_* \rangle \nonumber\\&\qquad \qquad \qquad \qquad \qquad \qquad \qquad-\lambda \langle F,\hat{S}-S_* \rangle -\lambda\langle \mathcal{P}_\Omega B ,\hat{S}-S_*\rangle \nonumber\\
&= G(L_*,S_*) + \langle X,\mathcal{P}_\mathcal{R}^\perp(\hat{L}-L_*)\rangle - \langle W,\hat{L}-L_*\rangle+\lambda \langle Y,\mathcal{P}_\Omega^\perp (\hat{S}-S_*) \rangle \nonumber\\&\qquad \qquad \qquad \qquad \qquad \qquad \quad-\lambda \langle F,\hat{S}-S_* \rangle -\lambda\langle \mathcal{P}_\Omega B ,\hat{S}-S_*\rangle \nonumber\\
&\labelrel \geq {myeq:3}  G(L_*,S_*) + \Vert \mathcal{P}_\mathcal{R}^\perp(\hat{L}-L_*)\Vert_{*} - \Vert  W\Vert \Vert \mathcal{P}_\mathcal{R}^\perp(\hat{L}-L_*)\Vert_{*}\nonumber\\&\qquad \qquad \qquad +\lambda \Vert \mathcal{P}_\Omega^\perp (\hat{S}-S_*)\Vert_{1} -\lambda \Vert  F\Vert_\infty \Vert \mathcal{P}_\Omega^{\perp}(\hat{S}-S_*) \Vert_{1} \nonumber\\&\qquad \qquad \qquad \qquad \qquad \qquad \qquad \qquad \quad -\lambda \Vert \mathcal{P}_\Omega B \Vert_{F} \Vert \mathcal{P}_\Omega (\hat{S}-S_*) \Vert_{F} \nonumber\\
&=   G(L_*,S_*) + \left( 1 - \Vert  W\Vert\right)  \Vert \mathcal{P}_\mathcal{R}^\perp(\hat{L}-L_*)\Vert_{*}\nonumber\\&\qquad \qquad \qquad  \qquad \qquad+\lambda \left( 1 - \Vert  F\Vert_\infty\right)  \Vert \mathcal{P}_\Omega^{\perp}(\hat{S}-S_*) \Vert_{1} \nonumber\\&\qquad \qquad \qquad  \qquad \qquad \qquad \qquad \qquad  -\lambda \Vert \mathcal{P}_\Omega B \Vert_{F} \Vert \mathcal{P}_\Omega (\hat{S}-S_*) \Vert_{F},\nonumber
\end{align}
where the equality in \eqref{myeq:11} follows by condition \ref{itm:first} of Lemma \ref{l2}, the equality in \eqref{myeq:1} follows by fact that both $(L_*,S_*)$ and $(\hat{L},\hat{S})$ satisfies the feasibility constraint, thus, $\langle UU^T+W,\hat{L}-L_*+\hat{S}-S_* \rangle =0$, and we make use of condition \ref{itm:snd} of Lemma \ref{l2}, that is, $\mathcal{P}_\mathcal{R}(X)=UU^T$, and $\mathcal{P}_\Omega (Y) = \sgn(C\circ S_*)$. The equality in \eqref{myeq:2} follows by $\mathcal{P}_\mathcal{R}(X)+\mathcal{P}_\mathcal{R}^\perp(X) = X$, and $\mathcal{P}_\Omega(Y)+\mathcal{P}_\Omega^\perp(Y) = Y$. The inequality in \eqref{myeq:3} follows from the fact that the dual norm of spectral norm is the nuclear norm, $ \langle W, \mathcal{P}_\mathcal{R}^\perp (\hat{L}-L_*) \rangle \leq  \Vert W\Vert \Vert \mathcal{P}_\mathcal{R}^\perp (\hat{L}-L_*)\Vert_{*}$, the details can be found in Lemma 3.2 in \cite{candes2009exact}, and the dual norm of the infinity norm is the $\ell_1$-norm \cite{yang1991generalized}, $ \langle F, \mathcal{P}_\Omega^\perp (\hat{S}-S_*) \rangle \leq  \Vert F\Vert_\infty \Vert \mathcal{P}_\Omega^\perp (\hat{S}-S_*)\Vert_{1}$. We have chosen $X$ and $Y$, such that $\langle X,\mathcal{P}_\mathcal{R}^\perp(\hat{L}-L_*) \rangle = \Vert \mathcal{P}_{\mathcal{R}^\perp}(\hat{L}-L_*)\Vert_{*}$ and $\langle Y ,\mathcal{P}_\Omega^\perp(\hat{S}-S_*) \rangle = \Vert \mathcal{P}_{\Omega}^\perp(\hat{S}-S_*)\Vert_{1}$.
Thus, with $\Vert W\Vert \leq \frac{\alpha}{2}$, $\Vert F\Vert_{\infty}<\frac{1}{2}$ and $\Vert \mathcal{P}_\Omega B\Vert_F \leq \frac{1}{4}$, for $0.0021 <\alpha<0.0914 $, we have 
\begin{align}\label{Done}
G(\hat{L},\hat{S})&\geq G(L_*,S_*) + \left( 1 - \frac{\alpha}{2}\right)  \Vert \mathcal{P}_\mathcal{R}^\perp(\hat{L}-L_*)\Vert_{*}+\frac{\lambda}{2}  \Vert \mathcal{P}_\Omega^{\perp}(\hat{S}-S_*) \Vert_{1}\nonumber\\&\qquad \qquad \qquad \qquad \qquad \qquad \qquad \qquad \qquad - \frac{\lambda}{4} \Vert \mathcal{P}_\Omega (\hat{S} - S_*)\Vert_F.
\end{align}
It follows from $\hat{S}-S_*=\mathcal{P}_\mathcal{R}(\hat{S}-S_*)+\mathcal{P}_\mathcal{R}^\perp (\hat{S}-S_*)$ that
\begin{align*}
\Vert \mathcal{P}_\Omega (\hat{S} -S_*)\Vert_F &\leq \Vert \mathcal{P}_\Omega \mathcal{P}_\mathcal{R} (\hat{S} -S_*)\Vert_F +\Vert \mathcal{P}_\Omega \mathcal{P}_\mathcal{R}^\perp (\hat{S}-S_*)\Vert_F\\
&\leq \frac{1}{2} \Vert \hat{S}-S_* \Vert_F +\Vert \mathcal{P}_\mathcal{R}^\perp (\hat{S}-S_*)\Vert_F\\
&\leq \frac{1}{2} \Vert \mathcal{P}_\Omega (\hat{S}-S_*) \Vert_F + \frac{1}{2} \Vert \mathcal{P}_\Omega^\perp (\hat{S}-S_*) \Vert_F +\Vert \mathcal{P}_\mathcal{R}^\perp (\hat{S}-S_*)\Vert_F.
\end{align*}
This implies that 
\begin{align*}
\frac{1}{2}\Vert \mathcal{P}_\Omega (\hat{S} -S_*)\Vert_F &\leq \frac{1}{2} \Vert \mathcal{P}_\Omega^\perp (\hat{S}-S_*) \Vert_F +\Vert \mathcal{P}_\mathcal{R}^\perp (\hat{S}-S_*)\Vert_F.
\end{align*}

Hence, it follows that 

\begin{align}
\frac{\lambda}{4} \Vert \mathcal{P}_\Omega (\hat{S}-S_*) \Vert_F  &\leq \frac{\lambda}{4} \Vert \mathcal{P}_\Omega^\perp (\hat{S}-S_*) \Vert_F +\frac{\lambda}{2}\Vert \mathcal{P}_\mathcal{R}^\perp (\hat{S}-S_*)\Vert_F\nonumber\\
&\leq \frac{\lambda}{4} \Vert \mathcal{P}_\Omega^\perp (\hat{S}-S_*) \Vert_1 +\frac{\lambda}{2}\Vert \mathcal{P}_\mathcal{R}^\perp (\hat{S}-S_*)\Vert_*\nonumber\\
&= \frac{\lambda}{4} \Vert \mathcal{P}_\Omega^\perp (\hat{S}-S_*) \Vert_1 + \frac{\lambda}{2}\Vert \mathcal{P}_\mathcal{R}^\perp (\hat{L}-L_*)\Vert_*, \nonumber
\end{align}
where the last equality follows by the feasibility of $(\hat{L},\hat{S})$ and $(L_*,S_*)$, that is $\hat{L} +\hat{S}=L_*+S_*$.
Therefore, from \eqref{Done} we have
\begin{align*}
G(\hat{L},\hat{S})&\geq G(L_*,S_*) + \left( 1 - \frac{\alpha}{2}\right)  \Vert \mathcal{P}_\mathcal{R}^\perp(\hat{L}-L_*)\Vert_{*}+\frac{\lambda}{2} \Vert \mathcal{P}_\Omega^{\perp}(\hat{S}-S_*) \Vert_{1} \nonumber\\&\qquad \qquad \qquad \qquad \qquad \qquad- \frac{\lambda}{4} \Vert \mathcal{P}_\Omega^\perp (\hat{S} - S_*)\Vert_1-\frac{\lambda}{2}\Vert \mathcal{P}_\mathcal{R}^\perp (\hat{L}-L_*)\Vert_*\\
&= G(L_*,S_*) + \left( 1 - \frac{\alpha+\lambda}{2}\right)  \Vert \mathcal{P}_\mathcal{R}^\perp(\hat{L}-L_*)\Vert_{*}+\frac{\lambda}{2} \Vert \mathcal{P}_\Omega^{\perp}(\hat{S}-S_*) \Vert_{1},
\end{align*}
given $\Vert W\Vert \leq \frac{\alpha}{2} $, and $\Vert F\Vert_{\infty}<\frac{1}{2}$, and we have $\frac{\alpha +\lambda}{2}<1$ as $0.0021 <\alpha< 0.0914$ and $\lambda<\alpha$, one has $G(\hat{L},\hat{S})> G(L_*,S_*)$, for $(\hat{L}-L_*, \hat{S}-S_*)\neq (0,0)$. 

However, $\Vert \mathcal{P}_\mathcal{R}^\perp(\hat{L}-L_*)\Vert_{*} = \Vert \mathcal{P}_\Omega^{\perp}(\hat{S}-S_*) \Vert_{1} =0 $  only if $\mathcal{P}_\mathcal{R}^\perp (\hat{L}-L_*)=\mathcal{P}_\Omega^\perp(\hat{S}-S_*) = 0$ (i.e., $\hat{L}-L_*=\hat{S}-S_* =P \in \Omega \cap \mathcal{R}$) then the injectivity assumption (that $\mathcal{R}$ and $\Omega$  have a trivial intersection) forces $(\hat{L}-L_*,\hat{S}-S_*)=(P,P)=(0,0)$.

Consequently, any minimizer $(\hat{L},\hat{S})$ with $(\hat{L}-L_* ,\hat{S}-S_*)\neq (0,0)$ must satisfies $G(\hat{L},\hat{S})>G(L_*,S_*)$. Thus, $(L_*,S_*)$ is a unique minimizer to problem \eqref{eq:m}. %\Halmos
\qed

According to \textbf{ Lemma \ref{l2}}, for the exact recovery of problem \eqref{eq:m}, it is sufficient to find an appropriate $W$, for which:
\qquad \qquad
\begin{equation}\label{eq:c1}
\left\{
\begin{array}{l}
W\in \mathcal{R}^\perp,\\
\Vert W\Vert \leq \frac{\alpha}{2},~0.0021 <\alpha< 0.0914 ,\\
\Vert \mathcal{P}_{\Omega}(UU^T+W-\lambda \sgn(C\circ S_*))\Vert_{F}\leq \frac{\lambda}{4},\\
\Vert \mathcal{P}_{\Omega}^\perp(UU^T+W)\Vert_{\infty}\leq \frac{\lambda}{2}.
\end{array}
\right.\\ 
\end{equation}

In the following section we use the Golfing scheme to construct an approximation of the dual certificate in the setting of decomposing a matrix into its low-rank and sparse components.  

\subsection{Dual certification by the Golfing scheme and least squares}
The main idea is to construct $W$ such that it satisfies the conditions in \eqref{eq:c1}. Assume that entries of $S$ are sampled according to Bernoulli model with probability $p$. This means that all the matrices in $\Omega^{\perp} \sim Ber(1-p)$. Assume that all the matrices in $\Omega^\perp$ has the same distribution as $\Omega_1 \cup \Omega_2\cup...\cup \Omega_K$, where $\Omega_k$'s are drown independently with replacement from $Ber(q)$, $K=20\lceil \log N \rceil$; the parameter $q$ is found empirically. This can be described by Binomial model, $Bin(K,q)$, that is,
$$\Pr((i,j) \in \Omega) = \Pr(Bin(K,q)=0)=(1-q)^K.$$
Thus, the two model are equivalent if $p=(1-q)^K$.

The key idea is to decompose $W$ into $W^L$ (low-rank component) and $W^S$ (sparse component), that is, $W = W^L + W^S$. Then
\begin{align*}
UU^T +W &= UU^T + W^L +W^S\\
& = \mathcal{P}_\Omega (UU^T + W^L +W^S) +\mathcal{P}_\Omega^\perp  (UU^T + W^L +W^S)\\
& = \mathcal{P}_\Omega (UU^T + W^L) +\lambda \sgn(C\circ S_*) +\mathcal{P}_\Omega^\perp  (UU^T + W^L +W^S)\\
& =  \lambda  \left( \mathcal{P}_\Omega \left( \frac{UU^T + W^L}{\lambda}\right)  +\sgn(C\circ S_*) \right)\\&\qquad \qquad \qquad \qquad \qquad\qquad+ \lambda \mathcal{P}_\Omega^\perp  \left( \frac{UU^T + W^L +W^S}{\lambda}\right) ,
\end{align*}
where $\mathcal{P}_\Omega (W^S)=\lambda \sgn(C\circ S_*)$, since $L_*$ and $S_*$ are supported on $\mathcal{R}$ and $\Omega$, respectively.
We take
\begin{align*}
B &= \left( \frac{UU^T + W^L}{\lambda}\right), \textit{ and } \\
F &= \mathcal{P}_\Omega^\perp  \left( \frac{UU^T + W^L +W^S}{\lambda}\right),
\end{align*} 
and so $B$ and $F$ adhere to the conditions stated in Lemma \ref{l2}.  $W^L$ and $W^S$ adhering to \eqref{eq:c1} certify that problem \eqref{eq:m} perfectly recovers the low-rank matrix $L_*$ and the sparse matrix $S_*$ with high probability, that is 

\begin{equation}\label{eq:cc1}
\left\{
\begin{array}{l}
\Vert W^L +W^S\Vert \leq \frac{\alpha}{2},~0.0021 <\alpha< 0.0914,\\
\Vert \mathcal{P}_{\Omega}(UU^T+W^L)\Vert_{F}\leq \frac{\lambda}{4},\\
\Vert \mathcal{P}_{\Omega}^\perp(UU^T+W^L+W^S)\Vert_{\infty}\leq \frac{\lambda}{2}.
\end{array}
\right.\\ 
\end{equation}

We will use the Golfing scheme to construct $W^L$, and the least squares method to construct $W^S$. The Golfing scheme \cite{gross2010quantum} is a tool to construct an approximate dual certificate. In the Golfing scheme, an interim solution is improved, iteratively, until the final approximation of the dual certificate is obtained. 
$W^L$ is constructed as follows: 
\begin{align}\label{1}
Q_k &= Q_{k-1} + q^{-1} \mathcal{P}_{\Omega_k} \mathcal{P}_{\mathcal{R}}(UU^T-Q_{k-1}),~ Q_0=0,~ k=1,2,...,K\\
&\qquad \qquad W^L =  \mathcal{P}_{\mathcal{R}}^{\perp} Q_K \nonumber
\end{align}

According to the least square method \cite{candes2011robust}, $W^S$ is constructed as follows:
\begin{align*}
W^S = \lambda \mathcal{P}_\mathcal{R}^\perp \left( \mathcal{P}_\Omega - \mathcal{P}_\Omega\mathcal{P}_\mathcal{R}\mathcal{P}_\Omega \right) ^{-1} \sgn(C\circ S_*),
\end{align*}
using Neumann series \cite{candes2011robust}, $W^S$ can be written as follows:
\begin{align}\label{WS}
W^S &= \lambda \mathcal{P}_{\mathcal{R}}^\perp \sum_{k=0}^{K} \left(  \mathcal{P}_{\Omega}  \mathcal{P}_{\mathcal{R}}  \mathcal{P}_{\Omega}\right)^k \sgn(C\circ S_*).
\end{align}

We now declare and verify some sufficient conditions on the approximated dual certificate $W^L$ constructed by the Golfing scheme for the pair $(L_*,S_*)$ to be the unique optimal solution to \eqref{eq:m}. 
\begin{lemma}\label{l3}
	Assume that all the matrices in $\Omega\sim Ber(p)$, $\Vert \mathcal{P}_\Omega \mathcal{P}_\mathcal{R}\Vert \leq \frac{1}{2}$, i.e., $\mathcal{P}_\Omega \cap \mathcal{P}_\mathcal{R}=0$, and $K=20\lceil \log N \rceil$. Then, for $\lambda =\frac{\alpha}{\sqrt{N}}$, $0.0021 <\alpha< 0.0914 $, the dual matrix $W^L$ in \eqref{1} satisfies:
	
	\begin{enumerate}
		\item[a.] $\Vert W^L \Vert < \frac{\alpha}{4}$,
		\item[b.] $\Vert \mathcal{P}_{\Omega}(UU^T + W^L) \Vert_F< \frac{\lambda}{4}$,
		\item[c.] $\Vert \mathcal{P}_{\Omega}^\perp(UU^T + W^L) \Vert_\infty < \frac{\lambda}{4}$.
	\end{enumerate}
\end{lemma}

The proof technique of this Lemma follows closely that of Lemma 2.8 in \cite{candes2011robust}, but we have have used different bounds in our proof.

{\it Proof}
\textbf{Proof of a.}  
Let 
\begin{align}\label{Z}
Y_{k-1} = UU^T -\mathcal{P}_\mathcal{R}  Q_{k-1}, ~Y_{k-1} \in \mathcal{R},
\end{align}
then, $Q_K =\sum_{k=1}^{K}  q^{-1} \mathcal{P}_{\Omega_k} \mathcal{P}_\mathcal{R} Y_{k-1}$, and $\Vert W^L \Vert = \Vert \mathcal{P}_{\mathcal{R}}^\perp  Q_K \Vert $. Note that it has been shown in \cite{candes2011robust} that, for
\begin{align}\label{q}
q\geq c_0 \zeta^{-2} \mu_0 r \log N /N,~~0<\zeta<1/2,
\end{align}
with high probability $ \Vert Y_k\Vert_{\infty} \leq \frac{1}{2} \Vert Y_{k-1} \Vert_{\infty} $ and $ \Vert Y_k\Vert_{F} \leq \frac{1}{2} \Vert Y_{k-1} \Vert_{F} $ hold, where $c_0$ is absolute constant. 

From the definition of $Q_k$ in \eqref{1} and using \eqref{Z} we have,
\begin{align}\label{TQ}
\mathcal{P}_\mathcal{R}^\perp Q_K=\sum_{k=1}^{K} q^{-1} \mathcal{P}_\mathcal{R}^\perp \mathcal{P}_{\Omega_k} \mathcal{P}_\mathcal{R} Y_{k-1}.
\end{align} 

Thus, using the following inequalities
\begin{align}\label{G}
\Vert \mathcal{P}_\mathcal{R}^\perp \mathcal{P}_{\Omega_k} Y_{k-1} \Vert \leq \frac{1}{4 \sqrt{\rank(Y_{k-1})}} \Vert Y_{k-1}\Vert_F,~~ \Vert Y_k\Vert_F \leq \frac{1}{2}\Vert Y_{k-1}\Vert_F,
\end{align}
found in  \cite{gross2010quantum},  we get:
\begin{align}\label{SAEQ}
\Vert W^L \Vert &= \Vert \mathcal{P}_{\mathcal{R}}^\perp  Q_K \Vert\\
& \leq  \sum_{k=1}^{K}  q^{-1}\Vert \mathcal{P}_\mathcal{R}^\perp \mathcal{P}_{\Omega_k}  Y_{k-1}\Vert,~~ \textit{by} ~~\eqref{TQ},\nonumber\\
&\leq q^{-1}\sum_{k=1}^{K} \frac{1}{4 \sqrt{a}}\Vert Y_{k-1}\Vert_F, \textit{ by }\eqref{G}, \textit{where the $\rank$ of } Y_{k-1}\textit{ is }a,\nonumber\\
&\labelrel \leq {SAEQ:S11} \frac{q^{-1}}{4\sqrt{a}} \sum_{k=1}^{K} \left( \frac{1}{2}\right) ^{k-1}\Vert UU^T\Vert_F, \textit{using } \Vert Y_{k-1}\Vert_F \leq \left( \frac{1}{2}\right) ^{k-1}\Vert Y_0\Vert_F \nonumber,\\
&\leq \frac{q^{-1}}{\sqrt{a}} \Vert UU^T\Vert_F,~~ \textit{since } \sum_{k=1}^{K} \left( \frac{1}{2}\right) ^{k-1}\leq 4, \textit{ for large } K \nonumber, \\
&\labelrel \leq {SAEQ:S22} \frac{C_0 \zeta^2 N \sqrt{N}}{\mu_0 r \sqrt{a}\log N} \Vert UU^T\Vert_\infty, \textit{  by } \eqref{q},\textit{ and the fact that }\Vert .\Vert_F \leq \sqrt{N}\Vert .\Vert_{\infty}\nonumber,\\
&\labelrel = {SAEQ:S33}\frac{C_0 \zeta^2 N^{\frac{3}{2}}}{\mu_0 r\sqrt{a} \log N} \frac{1}{n}\nonumber\\
&< \zeta^2 \alpha,\textit{ by choosing }C_0< \frac{\alpha n \mu_0 r \sqrt{a} \log N}{N^{\frac{3}{2}}},\nonumber\\
&< \frac{\alpha}{4},\textit{ by choosing }0<\zeta< \frac{1}{2}\nonumber.
\end{align}
where \eqref{SAEQ:S11} followed by \eqref{G} and using $Y_0 =UU^T$. We used $C_0=\frac{1}{c_0}$ in \eqref{SAEQ:S22}. \eqref{SAEQ:S33} followed by the joint incoherence condition in \eqref{eq:3} and using $\Vert UU^T\Vert_\infty=\frac{1}{n}$. 

\textbf{Proof of b. } From the definition of $Q_k$ in \eqref{1}, it is clear that $Q_K$ is supported on $\Omega^\perp$, this means that $\mathcal{P}_\Omega Q_K=0$. Thus, by the definition of $W^L$ we have
\begin{align}\label{W1}
\mathcal{P}_\Omega \left( UU^T +W^L \right) =  \mathcal{P}_\Omega \left( UU^T +\mathcal{P}_\mathcal{R}^\perp Q_K \right).
\end{align}

We also know that $\mathcal{P}_{\Omega} \left( \mathcal{P}_\mathcal{R}^\perp Q_K+\mathcal{P}_\mathcal{R} Q_K \right) = \mathcal{P}_\Omega (Q_K)=0$, this means that $$\mathcal{P}_\Omega (\mathcal{P}_\mathcal{R}^\perp Q_K)=-\mathcal{P}_\Omega (\mathcal{P}_\mathcal{R}Q_K),$$ and thus, substituting this in \eqref{W1} and using \eqref{Z}, we get
\begin{align*}
\mathcal{P}_\Omega \left( UU^T +W^L \right) =  \mathcal{P}_\Omega \left( UU^T -\mathcal{P}_\mathcal{R} Q_K \right) = \mathcal{P}_\Omega (Y_K).
\end{align*}

Therefore,
\begin{align}\label{Pb1}
\Vert \mathcal{P}_\Omega \left( UU^T +W^L \right) \Vert_F &= \Vert \mathcal{P}_\Omega Y_{K}\Vert_F\\
&\leq \Vert Y_{K}\Vert_F, ~ \textit{since } \Vert \mathcal{P}_{\Omega} Y_K\Vert_F \leq \Vert Y_K\Vert_F,\nonumber\\
&\leq \left( \frac{1}{2}\right) ^{K} \Vert UU^T\Vert_F,~ \Vert Y_K\Vert_F \leq \left( \frac{1}{2}\right) ^K\Vert Y_0 \Vert_F,~ Y_0 =UU^T,\nonumber\\
&\leq \left( \frac{1}{2}\right)^{K}  \sqrt{N}\Vert UU^T\Vert,\textit{ using the fact that } \Vert. \Vert_{F} \leq \sqrt{N}\Vert .\Vert,\nonumber\\
&= \left( \frac{1}{2}\right)^{K}  \sqrt{N},\textit{ as }\Vert UU^T\Vert =1,\nonumber\\
&= \frac{N}{\alpha ~2^K}\frac{\alpha}{\sqrt{N}}
< \frac{\lambda}{8}<\frac{\lambda}{4}, \textit{ where } \lambda=\frac{\alpha}{\sqrt{N}},\label{Pb2}
\end{align}
by choosing $K$ large enough such that $\frac{N}{\alpha~ 2^K} < \frac{1}{4}$, e.g., $K=20 \log N$.

\textbf{Proof of c. } We have 
\begin{align}
UU^T +W^L &= UU^T + \mathcal{P}_\mathcal{R}^{\perp}Q_K,~~ by ~~ \eqref{1},\nonumber\\ 
&= UU^T -\mathcal{P}_\mathcal{R} Q_K+ Q_K\nonumber \\
&= Y_K + Q_K. \label{QK}
\end{align}
Thus,
\begin{align}\label{SQ}
\Vert \mathcal{P}_{\Omega}^\perp  \left( UU^T +W^L\right)  \Vert_{\infty}  &= \Vert  \mathcal{P}_{\Omega}^\perp  \left( Y_K + Q_K\right)  \Vert_{\infty}, ~ \textit{by } \eqref{QK}, \\
&\leq \Vert Y_K\Vert_{\infty} +\Vert Q_K\Vert_{\infty} \nonumber\\
&\leq \Vert Y_K\Vert_{F} +\Vert Q_K\Vert_{\infty},~ \Vert Y_K\Vert_{\infty} \leq \Vert Y_K\Vert_{F}\nonumber,\\
&\leq \frac{\lambda}{8} +q^{-1}\sum_{k=1}^{K}\Vert \mathcal{P}_{\Omega} \mathcal{P}_\mathcal{R}Y_{k-1} \Vert_{\infty}\nonumber\\
&\leq \frac{\lambda}{8} +q^{-1}\frac{1}{2}\sum_{k=1}^{K}\Vert Y_{k-1} \Vert_{\infty},~\textit{as }\Vert \mathcal{P}_\Omega \mathcal{P}_\mathcal{R}  \Vert \leq\frac{1}{2},\nonumber\\
&\labelrel\leq {SQ:SA1} \frac{\lambda}{8} +q^{-1}\frac{1}{2}\sum_{k=1}^{K} (\frac{1}{2})^{k-1}\Vert UU^T \Vert_{\infty} \nonumber \\
&\labelrel={SQ:SA2} \frac{\lambda}{8} +2 q^{-1} \Vert UU^T \Vert_{\infty},~ \sum_{k=1}^{K} \left( \frac{1}{2}\right) ^{k-1} \leq 4,\textit{ for large }K, \nonumber\\
&\labelrel \leq {SQ:SA3} \frac{\lambda}{8} +2 \frac{C_0 \zeta ^2 N}{\mu_0 r \log N}  \frac{1}{n} \nonumber \\
& = \frac{\lambda}{8}  + 4  \frac{C_0 \frac{\zeta ^2}{2} N}{n\mu_0 r\log N} \left( \frac{\sqrt{N}}{\alpha} \right) \left( \frac{\alpha}{\sqrt{N}} \right) \nonumber \\
&< \frac{\lambda}{8} +\frac{\lambda}{8},\textit{ by choosing }0<\zeta< \frac{1}{2}, \textit{ thus, } \frac{\zeta ^2}{2} \frac{\alpha}{\sqrt{N}} <\frac{\lambda}{8},\nonumber
\end{align}
and by choosing $C_0< {\frac{ \alpha n \mu_0 r\log^2 N}{4N^{3/2} }}$. In the above derivation, the third inequality follows by the proof of b, i.e., the first inequality in \eqref{Pb1} and \eqref{Pb2}, along with \eqref{1} and \eqref{Z}. The ineuality in \eqref{SQ:SA1} followed by using $ \Vert Y_{k-1} \Vert_{\infty} \leq \left(\frac{1}{2} \right)^{k-1} \Vert Y_0 \Vert$, in \eqref{SQ:SA2} we used $Y_0=UU^T$, and the inequality \eqref{SQ:SA3} followed by \eqref{eq:3}, $\Vert UU^T\Vert_\infty=\frac{1}{n}, ~q \textit{ from } \eqref{q},~ C_0=\frac{1}{c_0}$.

 %\Halmos
\qed
We now declare and verify some sufficient conditions on the approximated dual certificate $W^S$ in Lemma \ref{l4}.
The following Lemma is somewhat similar to Lemma 2.9 in \cite{candes2011robust}, however, we have used the Bernstein's inequality, that provided a tighter bound than the one used in \cite{candes2011robust}.
\begin{lemma}[\cite{candes2011robust}]\label{l4}
	Assume that $S_*$ is supported on $\Omega$, and $\Vert \mathcal{P}_\Omega \mathcal{P}_\mathcal{R}  \Vert<\gamma$, $\gamma$ very small absolute number, with high probability. Then for $\lambda=\frac{\alpha}{\sqrt{N}}$, $0.0021 <\alpha< 0.0914$, the dual matrix $W^S$ in \eqref{WS} satisfies:
	
	\begin{enumerate}
		\item[a.] $\Vert W^S \Vert < \frac{\alpha}{4}$,
		\item[b.] $\Vert \mathcal{P}_{\Omega}^\perp( W^S) \Vert_\infty < \frac{\lambda}{4}$.
	\end{enumerate}
	
\end{lemma}

{\it Proof}
We consider the random variable $\delta_{ij} = \sgn ((C \circ S_{*})_{ij})$, such that:

\[
\delta_{ij}=
\begin{cases}
1, ~~ \textit{w.p } p,\\
0, ~~ \textit{w.p } 1-p.
\end{cases}	
\]
\newpage

\textbf{Proof of a.}

$W^S$ can be separated into two terms, using $k=0$ and $k\geq 1$ in \eqref{WS}:
\begin{align}
W^S	&=\lambda\mathcal{P}_\mathcal{R}^\perp \sgn(C \circ S_*)+  \lambda\mathcal{P}_\mathcal{R}^\perp \sum_{k= 1}^{K} \left( \mathcal{P}_{\Omega} \mathcal{P}_\mathcal{R} \mathcal{P}_\Omega \right)^k  \sgn(C \circ S_*)\label{WS1}.
\end{align}

The key idea of this prove is to find an upper bound of $\Vert W^S\Vert $. It follows that 

\begin{align}
\Vert \lambda\mathcal{P}_\mathcal{R}^\perp \sgn(C \circ S_*)\Vert \leq \lambda \Vert \sgn(C \circ S_*)\Vert \leq c\lambda\sqrt {Np}= c\alpha \sqrt{p}\leq \frac{\alpha}{8},\label{FBWS}
\end{align}
for small absolute constant $c$, where we have used the fact that in every dimension $N$, $\Vert \sgn(C \circ S_*)\Vert \leq c \sqrt{Np}$ \cite{vershynin2010introduction}, and we have $\lambda =\frac{\alpha}{\sqrt{N}}$.

We now consider the spectral norm of the second term of $W^S$ in \eqref{WS1}. Define $\mathcal{H}= \sum_{k= 1}^{K}  \left(  \mathcal{P}_{\Omega}  \mathcal{P}_{\mathcal{R}}  \mathcal{P}_{\Omega}\right)^k$ as an operator, then we can write 
$$\mathcal{P}_\mathcal{R}^\perp \sum_{k= 1}^{K}  \left( \mathcal{P}_{\Omega} \mathcal{P}_\mathcal{R} \mathcal{P}_\Omega \right)^k  \sgn(C \circ S_*),$$ as $ \mathcal{P}_\mathcal{R}^\perp \mathcal{H} \left(  \sgn(C \circ S_*)\right) $ and thus show that this term is bounded above by small absolute constant with high probability. 

Denote by $N$ the $\epsilon$-net of $\mathcal{S}^{N-1}$ of size at most $6^N$. According to \cite{vershynin2010introduction}, Lemma 5.4 (Computing the spectral norm on a net), we have

\begin{align*}
\Vert \mathcal{P}_\mathcal{R}^\perp \mathcal{H}(\sgn(C \circ S_*))\Vert &\leq \Vert \mathcal{H}(\sgn(C \circ S_*))\Vert \\
&= \sup _{x,y \in \mathcal{S}^{N-1}} \langle \mathcal{H} (yx^T), \sgn(C \circ S_*) \rangle\\
&\leq (1-2\epsilon)^{-1} \sup _{x,y \in N} \langle \mathcal{H} (yx^T), \sgn(C \circ S_*) \rangle\\
&= 4 \sup _{x,y \in N} \langle \mathcal{H} (yx^T), \sgn(C \circ S_*) \rangle,\textit{ using }\epsilon = \frac{3}{8}.
\end{align*}
Define the random variable $Z(x,y)=\langle \mathcal{H} (yx^T), \sgn(C \circ S_*) \rangle$, then by Matrix Bernstein's inequality \cite{tropp2015introduction}, for unit-normed vectors $x$ and $y$, that is, $\Vert x\Vert =\Vert y\Vert =1$, with zero mean, and variance $Var (Z(x,y))= \frac{1}{N^2} \sum_{N}\sum_{N} Z^2(x,y)= \frac{1}{N^2} \Vert Z(x,y)\Vert_{F}^2$, we have,

\begin{align*}
\Pr(\Vert Z(x,y)\Vert >a~\vert ~\Omega )&\leq  2N \exp \left(  \frac{-a^2/2}{\frac{1}{N^2} \Vert Z(x,y)\Vert_{F}^2+ a/3} \right),
\end{align*}
where $\Omega$ is the support of matrix $(\delta_{ij})$. Since $x$ and $y$ are unit-normed vectors, $\Vert yx^T\Vert_F=1$, $\Vert \mathcal{H} (yx^T)\Vert_{F}\leq \Vert \mathcal{H} \Vert $ thus,
\begin{align*}
\Pr(\Vert Z(x,y)\Vert >a ~| ~\Omega)&\leq 2N \exp \left(  \frac{-a^2/2}{\frac{1}{N^2} \Vert \mathcal{H}(yx^T)\Vert_{F}^2+ a /3} \right)\\
&\leq 2N \exp \left(  \frac{-a^2/2}{\frac{1}{N^2} \Vert \mathcal{H}\Vert^2+ a /3} \right).
\end{align*}

Therefore, we have 
\begin{align*}
\Pr(\lambda\Vert \mathcal{H}(C \circ \sgn(S_*))\Vert >a ~| ~\Omega)&\leq 2N \exp \left(  \frac{-(a/4\lambda)^2/2}{\frac{1}{N^2} \Vert \mathcal{H}\Vert^2+ (a/4\lambda) /3} \right)\\
&= 2N \exp \left(  \frac{-(a/\lambda)^2/32}{\frac{1}{N^2} \Vert \mathcal{H}\Vert^2+ a /12\lambda} \right).
\end{align*}

Assume that $\Vert \mathcal{P}_\Omega \mathcal{P}_\mathcal{R}\Vert \leq \gamma $ with high probability, for a very small absolute constant $\gamma$, we have 
\begin{align}\label{sigm}
\Vert\mathcal{H} \Vert &= \big \Vert \sum_{k= 1}^{K}   \left(  \mathcal{P}_{\Omega}  \mathcal{P}_{\mathcal{R}}  \mathcal{P}_{\Omega}\right)^k \big \Vert \\ &\leq \sum_{k= 1}^{K}  \big  \Vert \left(  \mathcal{P}_{\Omega}  \mathcal{P}_{\mathcal{R}}  \mathcal{P}_{\Omega}\right)^k \big \Vert \\
& \leq  \sum_{k= 1}^{K}  \gamma^{2k} = \frac{\gamma^2}{1-\gamma^2}\nonumber. 
\end{align}

Thus, unconditionally, 
\begin{align*}
\Pr(\lambda\Vert \mathcal{H}(\sgn(S_*))\Vert >a )&\leq 2N \exp \left( D  \right) \Pr(\Vert \mathcal{P}_\Omega \mathcal{P}_\mathcal{R}\Vert \leq \gamma)\\& \qquad \qquad \qquad \qquad \qquad + \Pr\left(  \Vert \mathcal{P}_\Omega \mathcal{P}_\mathcal{R}\Vert > \gamma \right) < \frac{\alpha}{8},
\end{align*}
where $\exp (D)$, with $D=\frac{-(a/\lambda)^2/32}{\frac{1}{N^2} \left( \frac{\gamma^2}{1-\gamma^2}\right) ^2+ a /12\lambda}$, is very small number, $\lambda =\frac{\alpha}{\sqrt{N}}$, and we put $a=\frac{\alpha}{8}$. This together with the bound of the first term of $\Vert W^S\Vert$ in \eqref{FBWS} completes the proof.  

\textbf{Proof b.} We know that $\mathcal{P}_\Omega W^S+\mathcal{P}_\Omega^\perp W^S=W^S$. Recalling $W^S$, we have

\begin{align}
\mathcal{P}_\Omega^\perp W^S &=W^S-\mathcal{P}_\Omega W^S\nonumber \\
&=\lambda \mathcal{P}_\mathcal{R}^\perp (\mathcal{I}-\mathcal{P}_\Omega) \sum_{k= 0}^{K} \left(  \mathcal{P}_{\Omega}  \mathcal{P}_{\mathcal{R}}  \mathcal{P}_{\Omega}\right)^k \sgn(C \circ S_*), \textit{ using } \eqref{WS}\nonumber,\\
&=\lambda\mathcal{P}_\Omega^\perp (\mathcal{I}-\mathcal{P}_\mathcal{R}) \sum_{k= 0}^{K}  \left(  \mathcal{P}_{\Omega}  \mathcal{P}_{\mathcal{R}}  \mathcal{P}_{\Omega}\right)^k \sgn(C \circ S_*)\nonumber\\
&= -\lambda \mathcal{P}_\Omega^\perp \mathcal{P}_{\mathcal{R}} \sum_{k=0}^{K} \left(  \mathcal{P}_{\Omega}  \mathcal{P}_{\mathcal{R}}  \mathcal{P}_{\Omega}\right)^k \sgn(C \circ S_*), \label{EQS}
\end{align}
where $\mathcal{I}$ is the identity operator and the last equality follows since $\sgn(C \circ S_*)$ is supported on $\mathcal{P}_\Omega$. The idea here is to express $ \Vert \mathcal{P}_\Omega^\perp W^S \Vert_{\infty}$ in the form of $\langle H, \sgn(C \circ S_*) \rangle$, then derive an upper bound on it, given $\Vert \mathcal{P}_\Omega \mathcal{P}_\mathcal{R} \Vert \leq \gamma$ (where $\gamma $ is a very small constant).

For any indices $(i,j)$ of $S$ $\in \Omega^\perp$, and noting that $\mathcal{P}_\Omega$ and $\mathcal{P}_\mathcal{R}$ are self ad-joint, thus 
\begin{align}
W_{ij}^S &= e_i^T W^S e_j=\langle e_i e_j^T, W^S \rangle \nonumber\\
& =\lambda\langle  e_i e_j^T, - \mathcal{P}_\Omega^\perp \mathcal{P}_{\mathcal{R}} \sum_{k=0}^{K} \left(  \mathcal{P}_{\Omega}  \mathcal{P}_{\mathcal{R}}  \mathcal{P}_{\Omega}\right)^k \sgn(C \circ S_*)  \rangle,~ \textit{using } \eqref{EQS},\nonumber\\
& =\lambda\langle  -\mathcal{P}_\Omega^\perp \mathcal{P}_{\mathcal{R}} \left( e_i e_j^T\right) ,  \sum_{k=0}^{K} \left(  \mathcal{P}_{\Omega}  \mathcal{P}_{\mathcal{R}}  \mathcal{P}_{\Omega}\right)^k \sgn(C \circ S_*)  \rangle \nonumber \\
& =\langle  - \sum_{k=0}^{K} \left(  \mathcal{P}_{\Omega}  \mathcal{P}_{\mathcal{R}}  \mathcal{P}_{\Omega}\right)^k\mathcal{P}_\Omega^\perp \mathcal{P}_{\mathcal{R}} \left( e_i e_j^T\right) ,  \lambda \sgn(C \circ S_*)  \rangle. \label{FFinal}
\end{align}

Define $Z(i,j) = - \sum_{k=0}^{K} \left(  \mathcal{P}_{\Omega}  \mathcal{P}_{\mathcal{R}}  \mathcal{P}_{\Omega}\right)^k\mathcal{P}_\Omega^\perp \mathcal{P}_{\mathcal{R}} \left( e_i e_j^T\right) $, thus using the union bound

\begin{align}
\Pr\left( \Vert \mathcal{P}_\Omega^\perp (W^S) \Vert > a\lambda~|~\Omega \right) &\leq \sum_{i,j} \Pr \left( \vert e_i^T W^S e_j \vert > a\lambda ~| ~ \Omega\right) \nonumber \\
&\leq N^2\Pr \left( \vert e_i^T W^S e_j \vert > a\lambda ~| ~ \Omega\right)\label{Final}.
\end{align}
Thus, using the matrix Bernstein's inequality, we have 

\begin{align*}
\Pr\left( \Vert \mathcal{P}_\Omega^\perp (W^S) \Vert_\infty > a\lambda~|~\Omega \right)&\leq 	\Pr\left(\sqrt{N} \Vert \mathcal{P}_\Omega^\perp (W^S) \Vert > a\lambda~|~\Omega \right),\textit{ as }\Vert .\Vert_{\infty} \leq \Vert .\Vert ,
\\
&\leq N^{5/2}\Pr \left( \vert e_i^T W^S e_j \vert > a\lambda ~| ~ \Omega\right),~ \textit{using } \eqref{Final},\\
& \leq  N^{5/2} \Pr \left( \vert \langle Z(i,j),\sgn(S_*) \rangle \vert >a ~|~\Omega\right), ~by\eqref{FFinal},\\
& \leq 2 N^{5/2 } \exp  \left( \frac{-a^2/2}{\frac{1}{N^2} \Vert Z(i,j) \Vert_F^2 + a/3}  \right),
\end{align*}
where the last inequality follows by the matrix Bernstein's inequality.  
Now for any indices $(i,j)$ of $S$ $\in \Omega^\perp$, assume that $\Vert \mathcal{P}_\Omega \mathcal{P}_\mathcal{R}\Vert \leq \gamma$, $\gamma$ small absolute number, then $\Vert \mathcal{P}_\Omega \mathcal{P} _\mathcal{R} \mathcal{P} _\Omega\Vert \leq \gamma^2$, thus we have 
\begin{align*}
\Vert Z(i,j)\Vert_{F}& =\sum_{k=0}^{K}\Vert  \left(  \mathcal{P}_{\Omega}  \mathcal{P}_{\mathcal{R}}  \mathcal{P}_{\Omega}\right)^k \mathcal{P}_\Omega^\perp \mathcal{P}_\mathcal{R}(e_ie_j^T) \Vert\\
&\leq \sum_{k=0}^{K}\Vert  \left(  \mathcal{P}_{\Omega}  \mathcal{P}_{\mathcal{R}}  \mathcal{P}_{\Omega}\right)^k \Vert \Vert  \mathcal{P}_\Omega^\perp \mathcal{P}_\mathcal{R}(e_ie_j^T) \Vert\\
&\leq \frac{\gamma^2}{1-\gamma^2} \Vert \mathcal{P}_\mathcal{R} (e_ie_j^T)\Vert_F,~ \textit{since }  \Vert  \mathcal{P}_\Omega^\perp \mathcal{P}_\mathcal{R}(e_ie_j^T) \Vert \leq  \Vert \mathcal{P}_\mathcal{R}(e_ie_j^T) \Vert,\\
&\leq \frac{\gamma^2}{1-\gamma^2} \sqrt {1-\Vert \mathcal{P}_\mathcal{R}^\perp (e_ie_j^T)\Vert^2_F}\\ 
&\leq \frac{\gamma^2}{1-\gamma^2} \sqrt {1-\Vert (I -\tilde{U}\tilde{U}^T)e_i\Vert^2\Vert (I -{\tilde{V}}\tilde{V}^T)e_j\Vert^2} \\
&= \frac{\gamma^2}{1-\gamma^2},~ \tilde{U}\tilde{U}^T=I,~ \tilde{V}\tilde{V}^T=I,
\end{align*}
where the third inequality follows from the fact that  $\Vert \mathcal{P}_\mathcal{R}(e_ie_j)\Vert_{F}^2 +\Vert \mathcal{P}_\mathcal{R}^\perp(e_ie_j)\Vert_{F}^2=1 $, the fourth inequality follows by the definition of the orthogonal complement projection onto $\mathcal{R}$ and using the fact that $, e_ie_j^T\textit{ has SVD } \tilde{U} \Sigma \tilde{V^T}$. Thus, unconditionally, 

\begin{align*}
\Pr\left( \Vert \mathcal{P}_\Omega^\perp (W^S) \Vert_\infty > a\lambda\right)
& \leq 2 N^{5/2 } \exp  \left( \frac{-a^2/2}{\frac{1}{N^2} \Vert Z(i,j) \Vert_F^2 + a/3}  \right) \Pr(\Vert \mathcal{P}_\Omega \mathcal{P}_\mathcal{R}\Vert \leq \gamma ) \\ & \qquad \qquad \qquad \qquad \qquad \qquad \qquad \quad + \Pr(\Vert \mathcal{P}_\Omega \mathcal{P}_\mathcal{R}\Vert > \gamma )\\
&\leq  2 N^{5/2 } \exp  \left( G  \right) \Pr(\Vert \mathcal{P}_\Omega \mathcal{P}_\mathcal{R}\Vert \leq \gamma ) + \Pr(\Vert \mathcal{P}_\Omega \mathcal{P}_\mathcal{R}\Vert > \gamma )\\
&< \frac{\lambda}{4},
\end{align*}	
where $\exp  \left( G \right)$, with $G=\frac{-a^2/2}{\frac{1}{N^2} \left( \frac{\gamma^2}{1-\gamma^2}\right)^2   + a/3} $, is a very small number,
 $\Pr\left( \Vert \mathcal{P}_\Omega \mathcal{P}_\mathcal{R}\Vert\leq \gamma  \right)$  with high probability and using $a=\frac{\lambda}{4}$, and $\lambda =\frac{\alpha}{\sqrt{N}}$. %\Halmos
\qed

The fact that $W^L$ and $W^S$ adhere to Lemma \ref{l3} and Lemma \ref{l4}, respectively, certifies that, with high probability, problem \eqref{eq:m} correctly recovers $L_*$ and $S_*$. 

\section{Numerical results} \label{Sec: Section 5}
In this section, we evaluate the performance of the proposed algorithm for solving the planted as well as maximum clique problem. All results are computed in Matlab 2019b, using a standard desktop computer with an Intel Core i7, 3.60GHz CPU, and 16 GB RAM. Here, we evaluate the performance of \textbf{Algorithm} \ref{a4} by applying it to identify the planted cliques in given graphs, to find maximum cliques in random graphs (where no cliques are planted), and finally to identify cliques for the real-world graphs. 

\subsection{Planted cliques}
Let $V^*$ denote the planted clique of size $n$. Let $M$ represent the adjacency matrix of the graph $G(V,E)$, $\vert V\vert =N$. We set $M_{ij} =1$ for $(i,j) \in V^* \times V^*$, $M_{ii}=1$ for all $i$; we add an edge $(i,j)$ with probability $p$ for all $(i,j) \in (V\times V)\backslash (V^* \times V^*) $ such that $M_{ji}=M_{ij}$. \textbf{Algorithm} \ref{a4} solves all problems tested with $p\in [\frac{1}{2}, 0.85]$ and achieves very similar accuracies for all $p$ values, see \cite{meunpub2020}. However, we report here the results obtained for $p=\frac{1}{2}$, for making a fair comparison with other algorithms in the literature. 

We have used $N=200,~500$ and $1000$. For a fixed value of $N$, we have used $n =10,20,...,N-10$.
Hence the number of problems considered for $N=200,~500$, and $1000$, are $19,~49,$ and $99$ respectively. Each of these problem is generated 15 times and hence the total number of test runs was 2505.

We have initialized \textbf{Algorithm} \ref{a4} with a randomly generated feasible $S$ of zeros with probability $p=0.75$ and ones with probability $1-p=0.25$. Then we initialize the feasible $L$ as $L=M-S$, $C=\frac{\epsilon}{(S+\epsilon)^2}$, and we set $l=1$. We would like to note here that the infeasible initialization such as $(L,S)=(0,0)$ equally produces similar final results. The parameters involved are $\delta$, and $\lambda$, where $\delta = 0.0001 $ is used as the tolerance for stopping the algorithm. We have used a constant $\lambda$ throughout our numerical testing. Our numerical investigations suggest that \textbf{Algorithm} \ref{a4} produces almost insensitive results for $\lambda = \frac{\alpha}{\sqrt{N}}$ for any $\alpha \in [ 0.0021 , 0.0914]$, see \cite{meunpub2020}. We have estimated the range, $[l,u]$, for $\alpha$ as follows. First we calculate three ranges $[l_i,u_i]$, $i=1,2,3$, corresponding to $c=0.25,0.5 \textit{ and } 0.75$, respectively in $m=\frac{1}{2}(N^2 – n^2)$, $n=cN$. We plot $\alpha$ ($=N/m$) against $N$ for each $c$ value and obtained $[l_i,u_i]$ for $\alpha$. We then take $l=\min{l_i}$ and $u=\max{u_i}$, $i=1,2,3$. We have used $\alpha = 0.054$ for all $(N,n)$ pairs for the results presented here. 
The suitable values of $\epsilon $ in \eqref{p1} lie in  $ [0.05,0.42]$, see \cite{meunpub2020}. For the results presented here we have used $\epsilon = 0.05$.  
We have used $\rho =\frac{1}{mean(M)}$, where $mean(M)$ is the mean value of entries of $M$. The regular version \eqref{e11}-\eqref{eq12} has been also implemented with these parameter values.

The final solution of ADMM algorithm for the regular model is denoted as $(L^1,S^1)$ while the final solution of the proposed model \eqref{eq9}-\eqref{eq91} is denoted as $(L^2,S^2)$.

We use the Frobenius norm to calculate the relative error 
$err L^i$ for each algorithm,
\begin{align}\label{E1}
errL^i=\frac{\Vert L^i-L_*\Vert _F}{\Vert L_*\Vert _F},~ i=1,2,
\end{align}
where $L_*$ corresponds to $V^*$, the planted clique.

We terminate \textbf{Algorithm} \ref{a4} when 
\begin{align}\label{SC1}
\Vert M - L_J - S_J \Vert_F \leq \delta,
\end{align}
holds. We have compared \textbf{Algorithm} \ref{a4} with the the densest subgraph algorithm (DSA) \cite{bombina2020convex,ames2015guaranteed}, for all the problems considered in this section.

\begin{figure}[H]
	\centering
	\subfloat[Subfigure 1 list of figures text][$N=200$]
	{
		\includegraphics[width=0.5\textwidth]{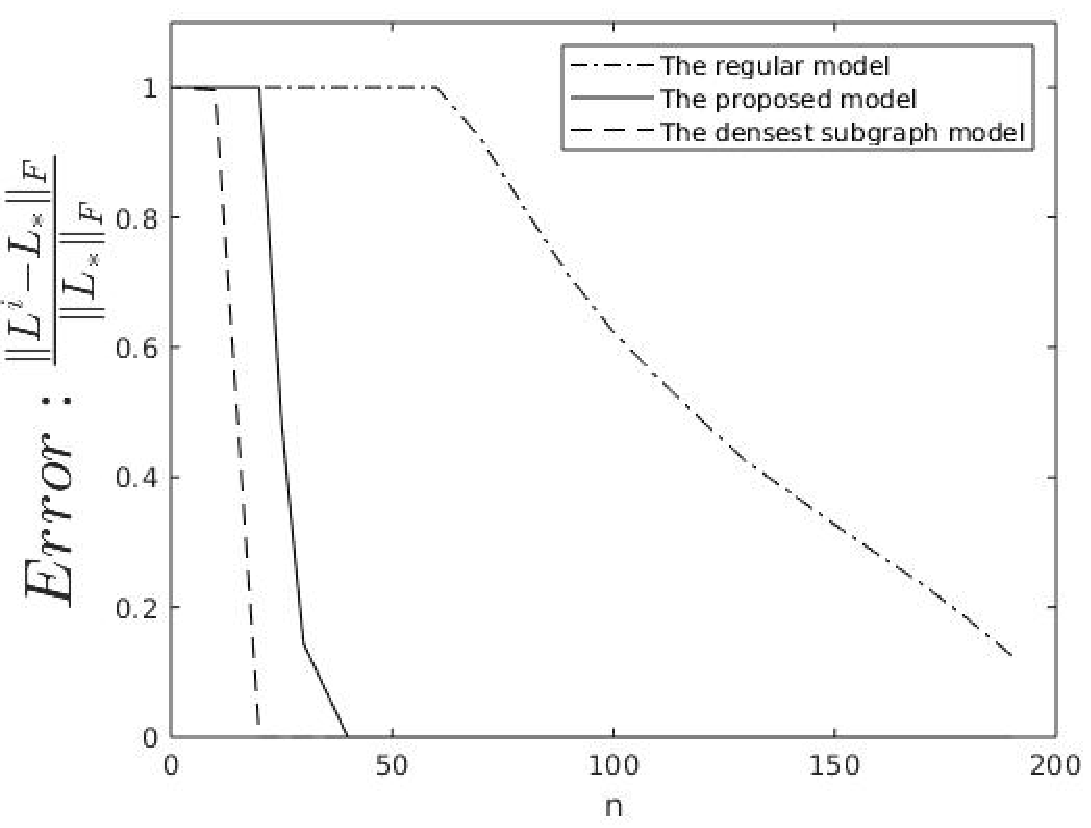}
		\label{fig:subfig1}
	}
	\subfloat[Subfigure 2 list of figures text][$N=500$]
	{
		\includegraphics[width=0.5\textwidth]{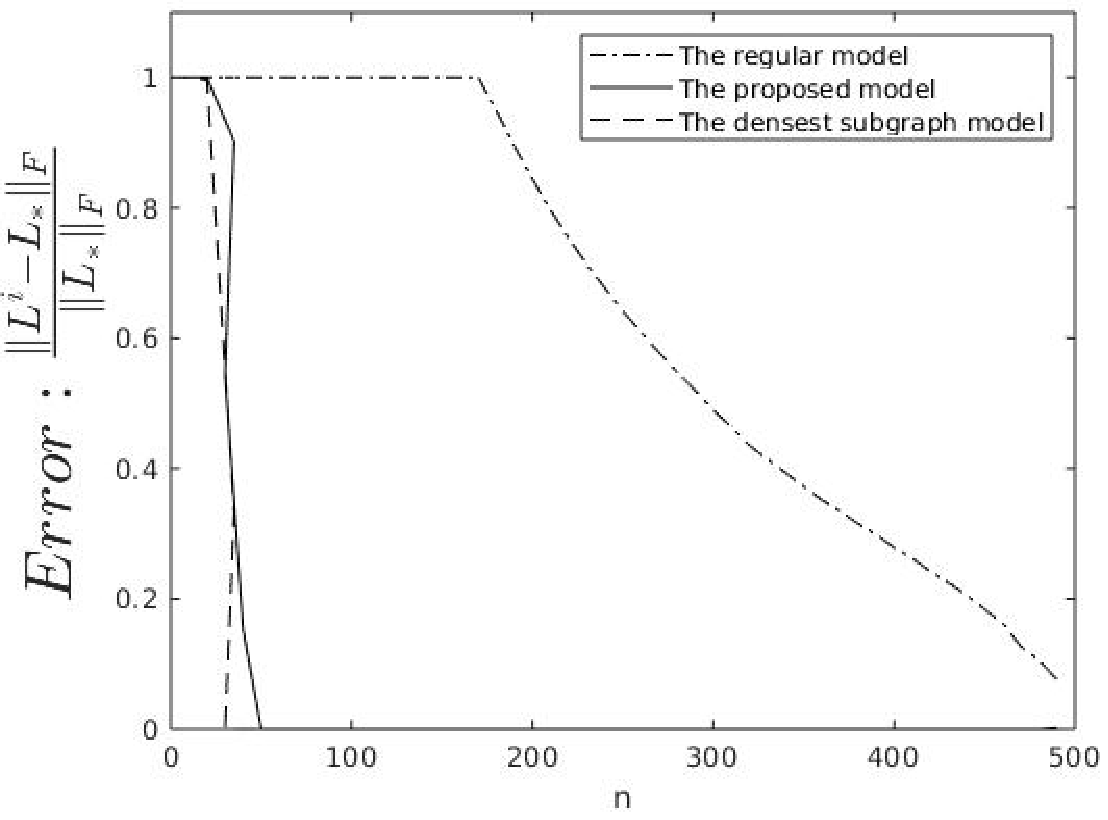}
		\label{fig:subfig2}
	}\\
	\subfloat[Subfigure 3 list of figures text][$N=1000$]
	{
		\includegraphics[width=0.5\textwidth]{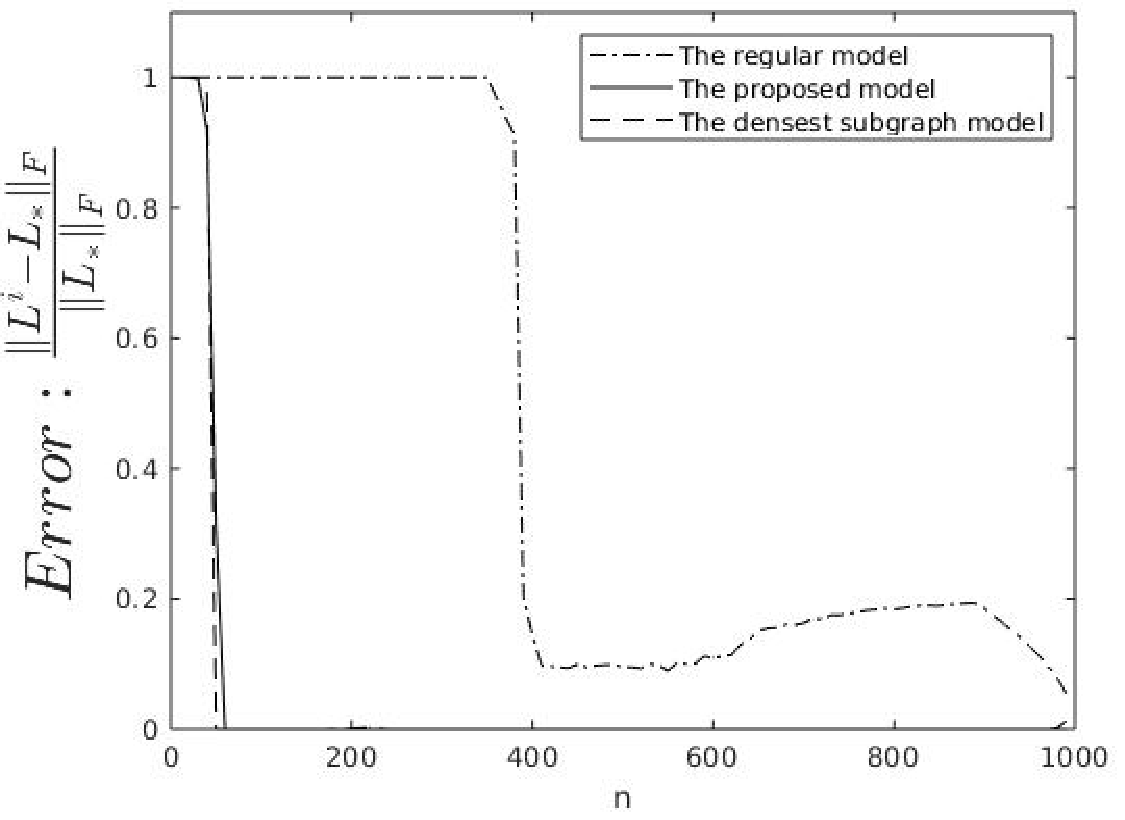}
		\label{fig:subfig3}
	}
	\caption{Average error for all problems.}
	\label{fig:fig1}
\end{figure}

 First we compare the average errors in Figure \ref{fig:fig1}, where the $y$-axis denotes the average of relative errors in \eqref{E1}; the average is taken over 15 runs on each problem. The value $n$ in the $x$-axis denotes the size of the planted clique. 
Figure \ref{fig:fig1} shows errors for the ADMM algorithm for the regular model are worse than the proposed model for all $(N,n)$ pairs. It also shows that error for the ADMM algorithm for the regular model does improve for higher values of $n$, i.e., for the problems that are easier to solve. On the other hand, our proposed model \eqref{eq9}-\eqref{eq91} achieves errors less than $10^{-8}$ for all $n\geq 30$, $n\geq 50$, and $n\geq 60$ for $N=200,~500$, and $1000$, respectively. However, the errors produced by DSA are about $10^{-5}$, it fails to produce less error than $10^{-5}$.

\begin{figure}[H]
	\centering
	\subfloat[Subfigure 1 list of figures text][$N=200$]
	{
		\includegraphics[width=0.5\textwidth]{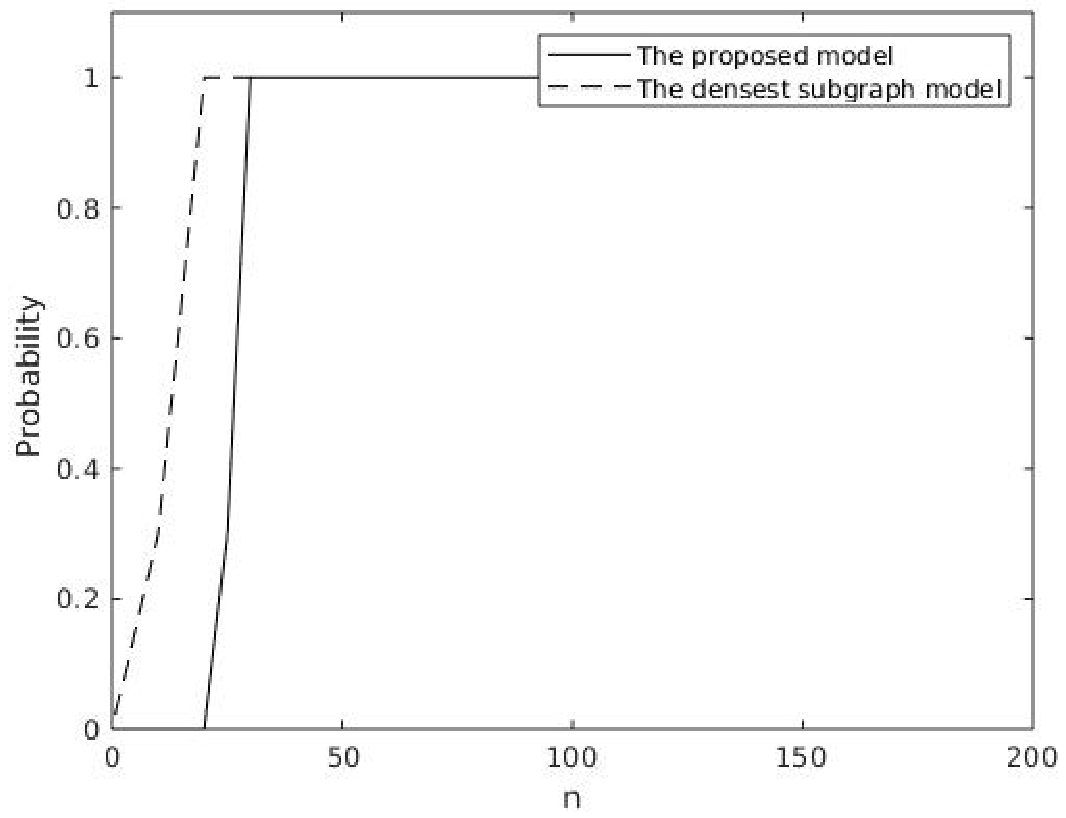}
		\label{fig:subfi1}
	}
	\subfloat[Subfigure 2 list of figures text][$N=500$]
	{
		\includegraphics[width=0.5\textwidth]{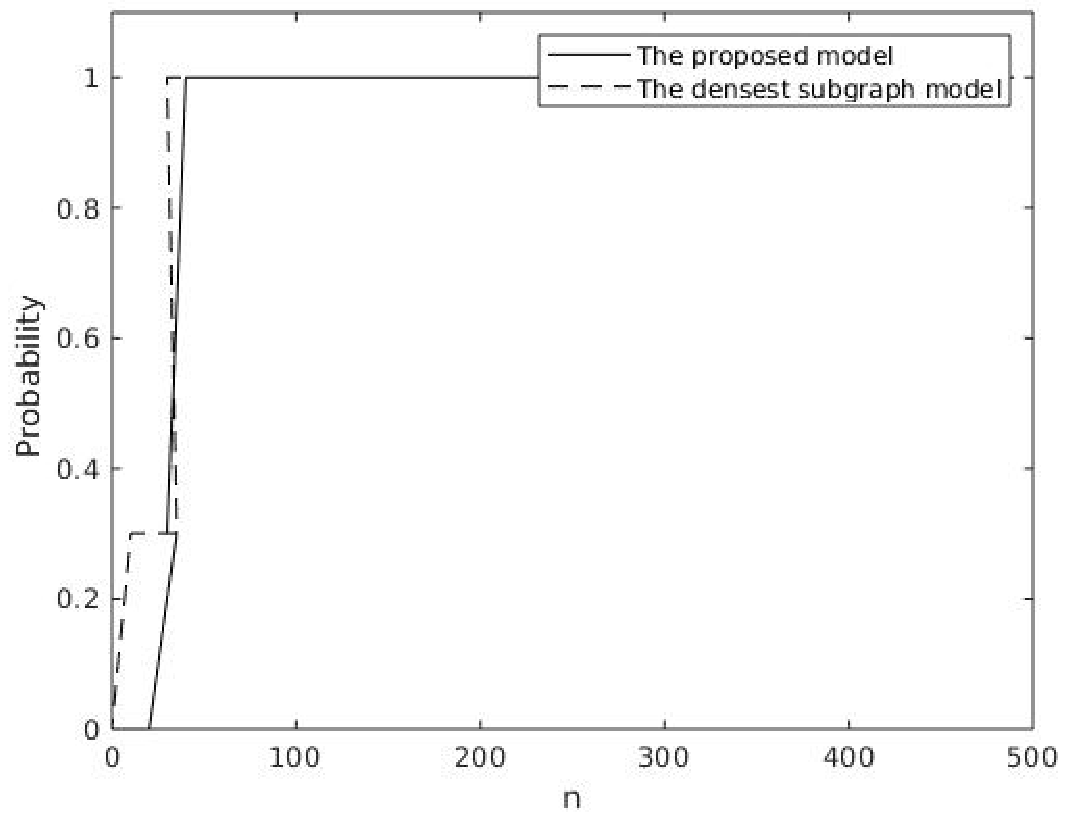}
		\label{fig:subfi2}
	}\\
	\subfloat[Subfigure 3 list of figures text][$N=1000$]
	{
		\includegraphics[width=0.5\textwidth]{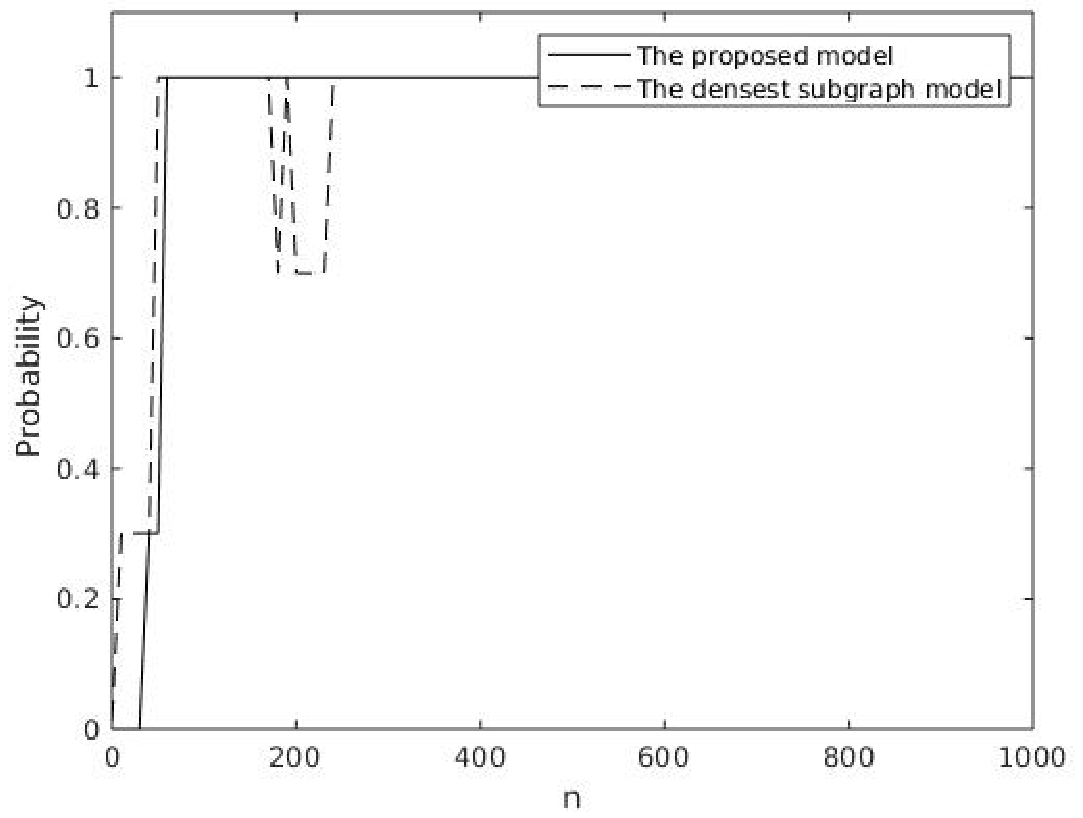}
		\label{fig:subfi3}
	}
	\caption{Average recovery probability for all problems.}
	\label{fig:fig2}
\end{figure}

Next, we compare our algorithm with the DSA using the probability of recovery for all problems corresponding to all $(N,n)$ pairs. This comparison has been summarized in Figure \ref{fig:fig2}. Here, by recovery we mean that the obtained solution has average error less than $10^{-8}$ for \textbf{Algorithm} \ref{a4} and about $10^{-5}$ for the DSA \cite{bombina2020convex,ames2015guaranteed}.  

Figure \ref{fig:fig2} shows that the probability equals one almost for all $(N,n)$ pairs using \textbf{Algorithm} \ref{a4}.
Figure \ref{fig:fig2} also shows that DSA has not provided perfect recovery for all the problems considered. For example, for $N=1000$, DSA has recovered some cliques of sizes around $n=200$ with probability less than 1 for a number of problems.  

\begin{figure}[H]
	\centering
	\subfloat[Subfigure 1 list of figures text][$N=200$]
	{
		\includegraphics[width=0.5\textwidth]{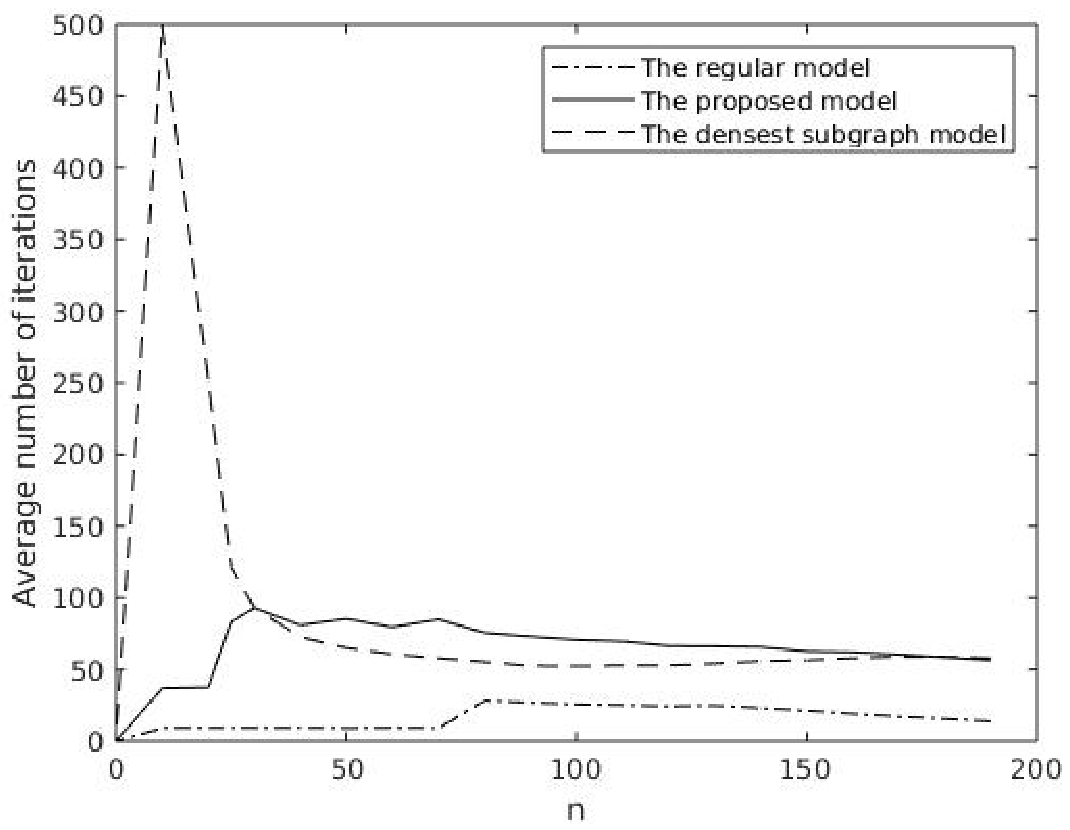}
		\label{fig:subfk1}
	}
	\subfloat[Subfigure 2 list of figures text][$N=500$]
	{
		\includegraphics[width=0.5\textwidth]{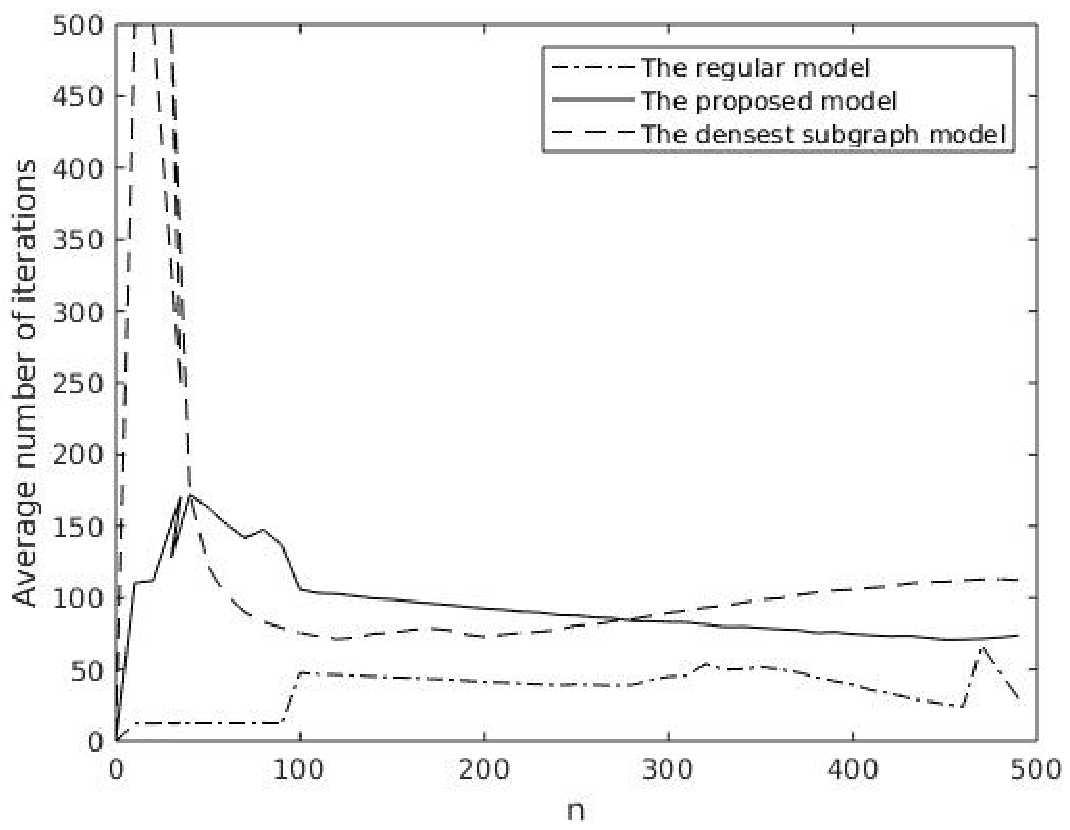}
		\label{fig:subfk2}
	}\\
	\subfloat[Subfigure 3 list of figures text][$N=1000$]
	{
		\includegraphics[width=0.5\textwidth]{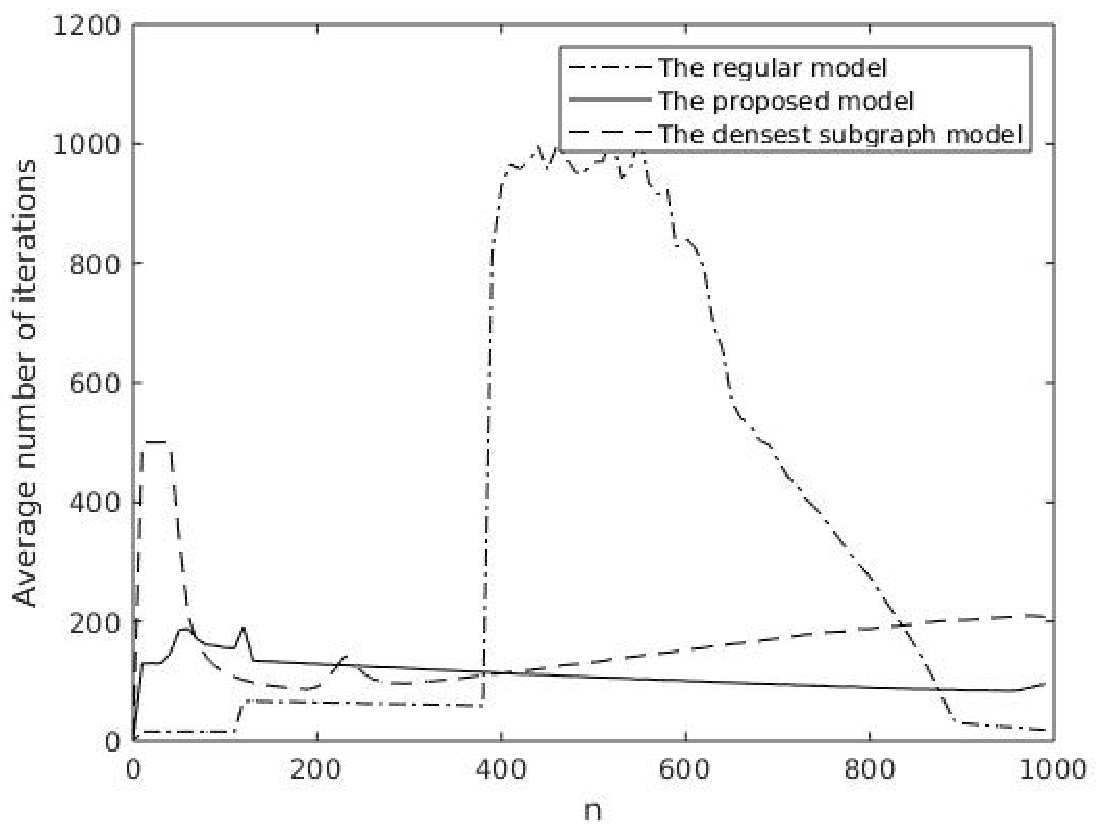}
		\label{fig:subfk3}
	}
	\caption{Average number of iterations for all problems.}
	\label{fig:fik2}
\end{figure}

In Figure \ref{fig:fik2} we present the average number of iterations needed by Algorithm \ref{a4} for producing average error of $10^{-8}$ and DSA for producing average error of $10^{-5}$.

Figure \ref{fig:fik2} demonstrates that both our proposed algorithm and DSA require comparable number of iterations to converge, except for high values of $(N,n)$ when DSA requires more iterations.  

To give an impression of the runtime needed by our algorithm and DSA, we have summarized the average of the total runtime and average runtime per iterations in the following figures.

\begin{figure}[H]
	\centering
	\subfloat[Subfigure 1 list of figures text][$N=200$]
	{
		\includegraphics[width=0.5\textwidth]{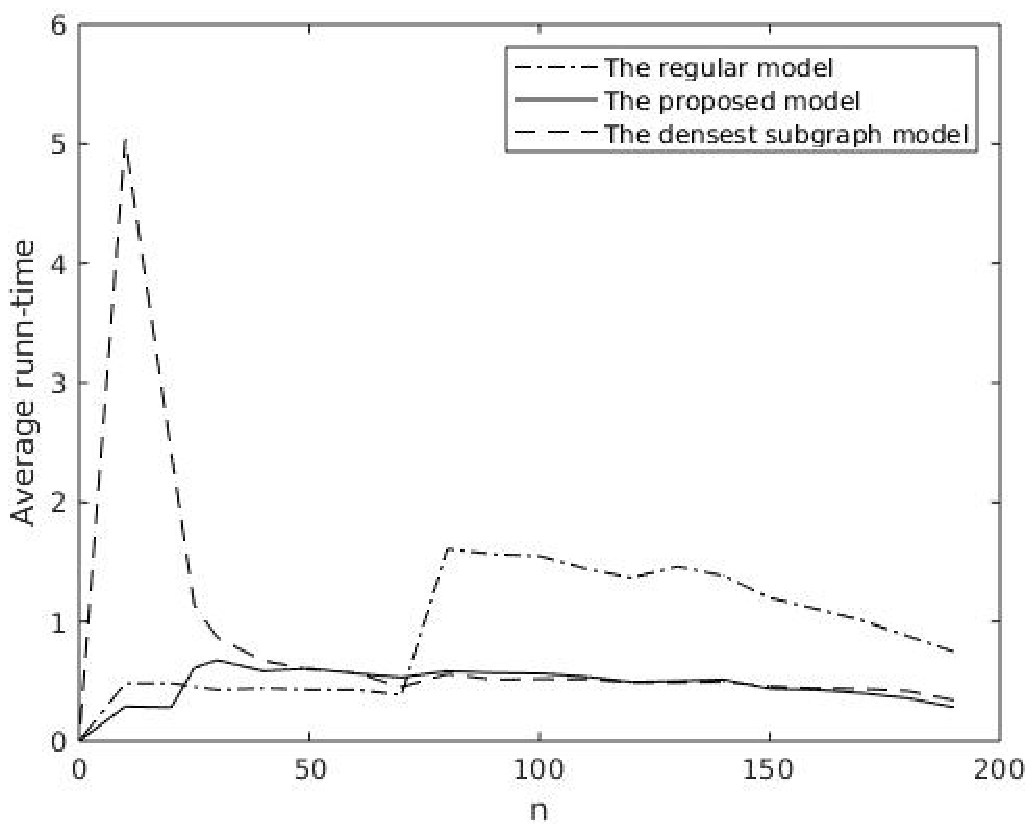}
		\label{fig:subf1}
	}
	\subfloat[Subfigure 2 list of figures text][$N=500$]
	{
		\includegraphics[width=0.5\textwidth]{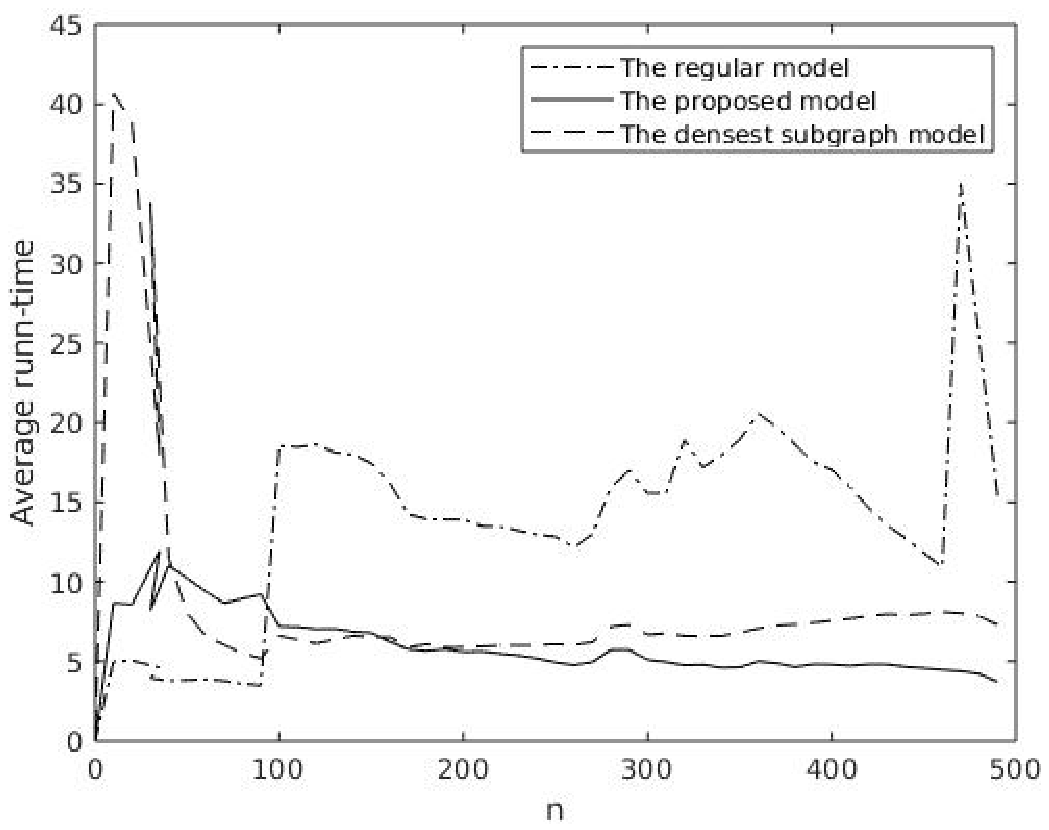}
		\label{fig:subf2}
	}\\
	\subfloat[Subfigure 3 list of figures text][$N=1000$]
	{
		\includegraphics[width=0.5\textwidth]{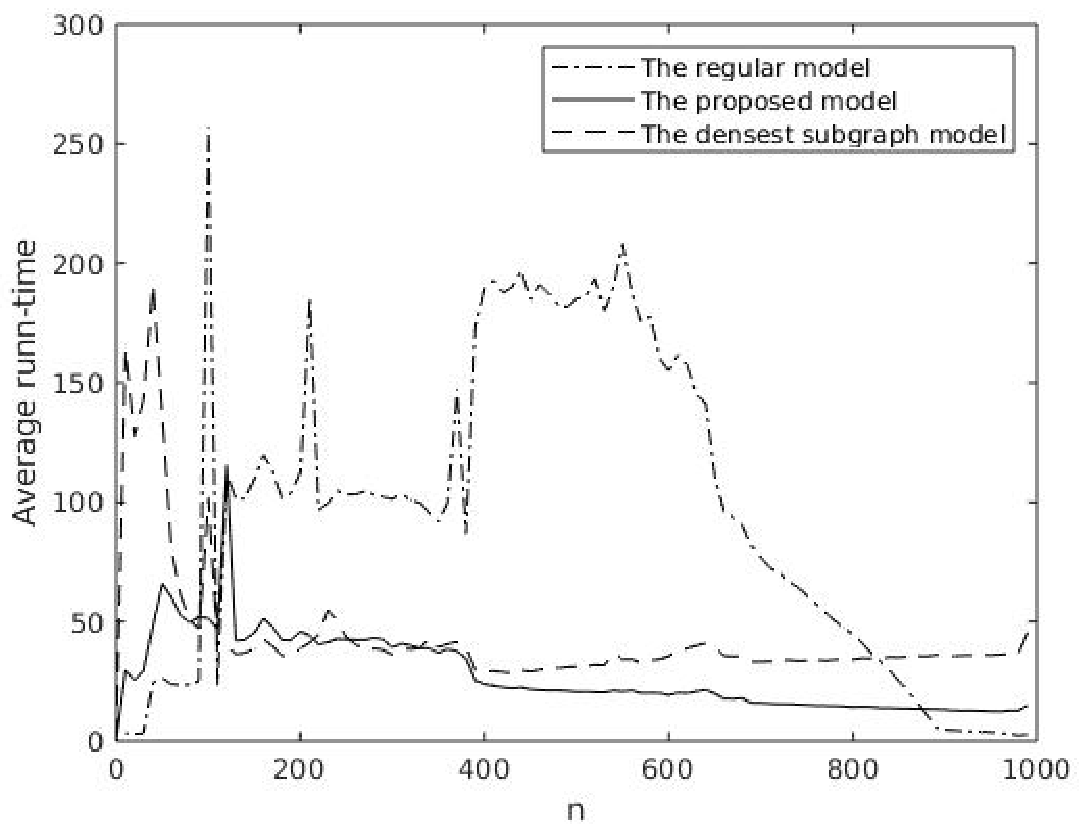}
		\label{fig:subf3}
	}
	\caption{Average runtime for all problems.}
	\label{fig:fi2}
\end{figure}
\begin{figure}[H]
	\centering
	\subfloat[Subfigure 1 list of figures text][$N=200$]
	{
		\includegraphics[width=0.5\textwidth]{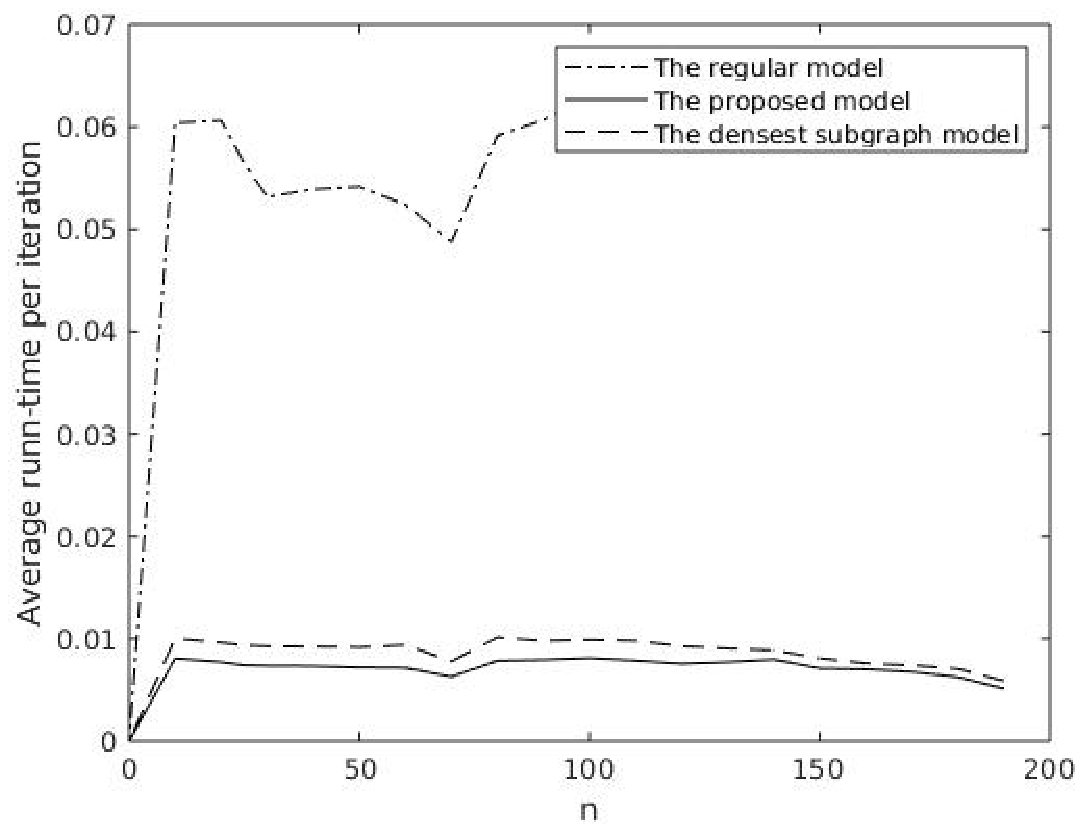}
		\label{fig:sub10}
	}
	\subfloat[Subfigure 2 list of figures text][$N=500$]
	{
		\includegraphics[width=0.5\textwidth]{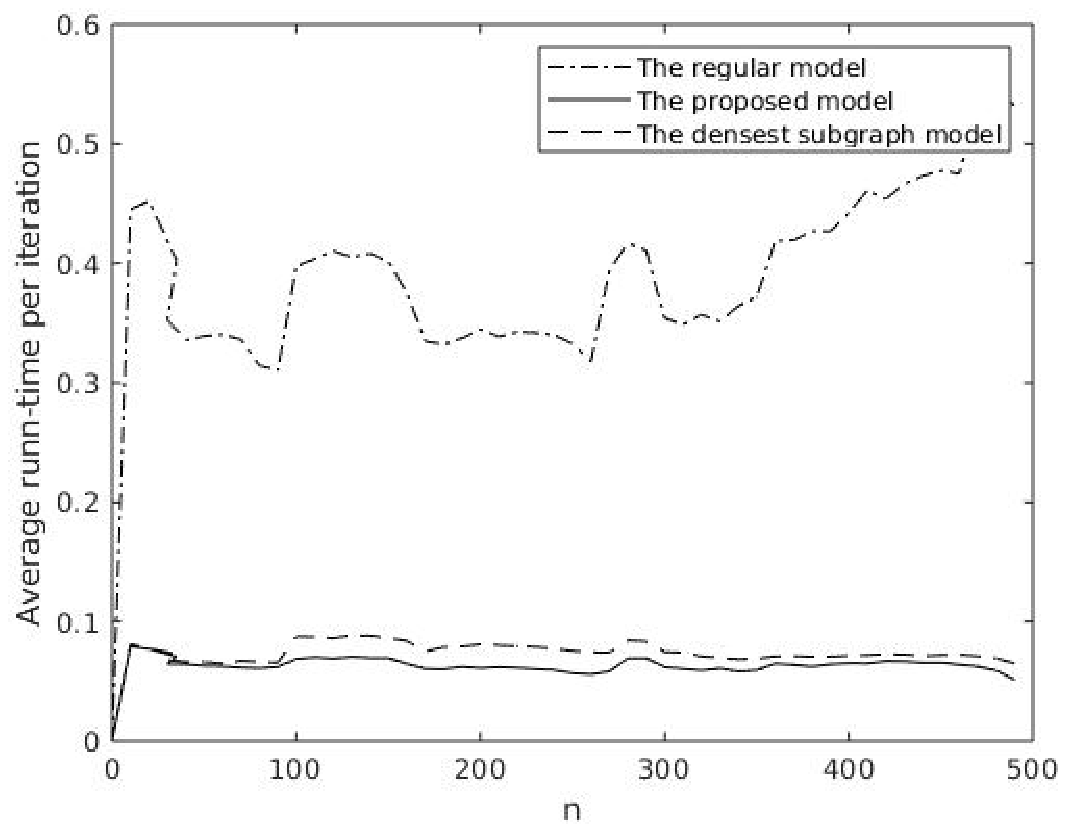}
		\label{fig:sub20}
	}\\
	\subfloat[Subfigure 3 list of figures text][$N=1000$]
	{
		\includegraphics[width=0.5\textwidth]{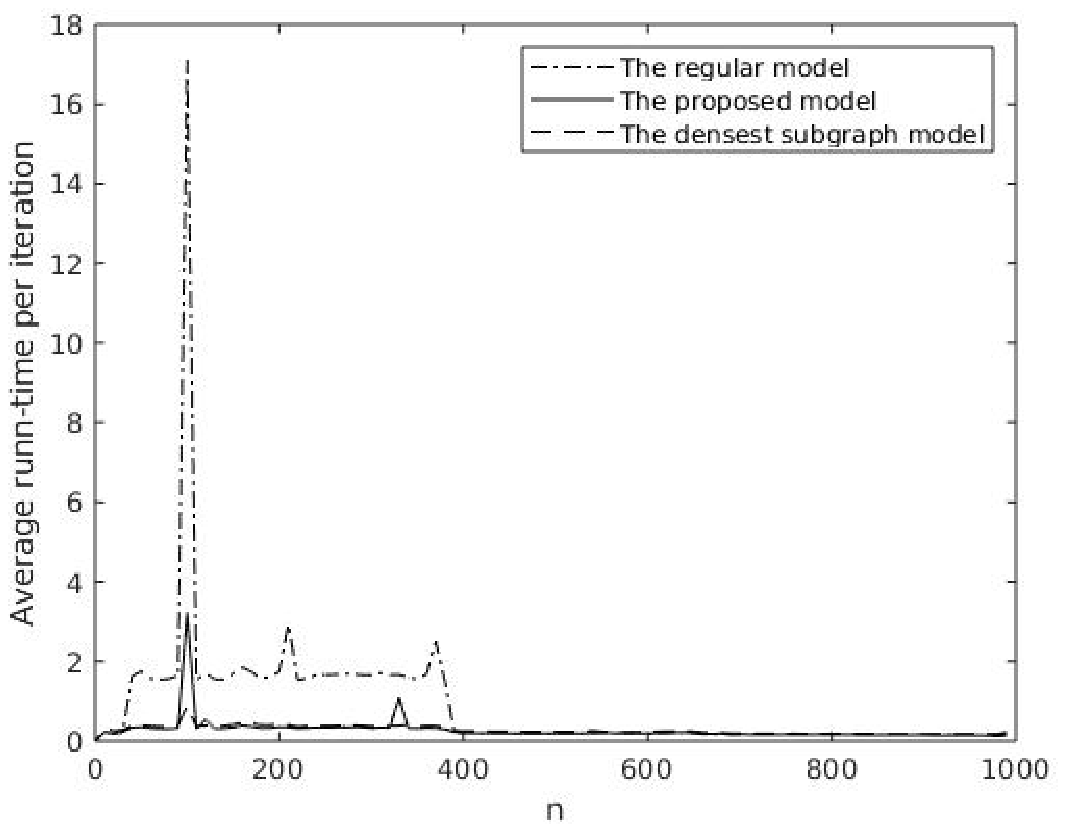}
		\label{fig:sub30}
	}
	\caption{Average runtime per iteration for all problems.}
	\label{fig:f2}
\end{figure}

Figure \ref{fig:fi2} shows the average runtime for all $(N,n)$ pairs for the problems considered, the averages are taken over 15 runs on each problem. Figure \ref{fig:fi2} shows that our algorithm performs better than DSA in finding the optimal solution. 
Figure \ref{fig:f2} shows the average runtime per iteration for the all the problems considered, where our algorithm performs slightly better.

\begin{figure}[H]
	\centering
	\subfloat[Subfigure 1 list of figures text][$N=200$]
	{
		\includegraphics[width=0.5\textwidth]{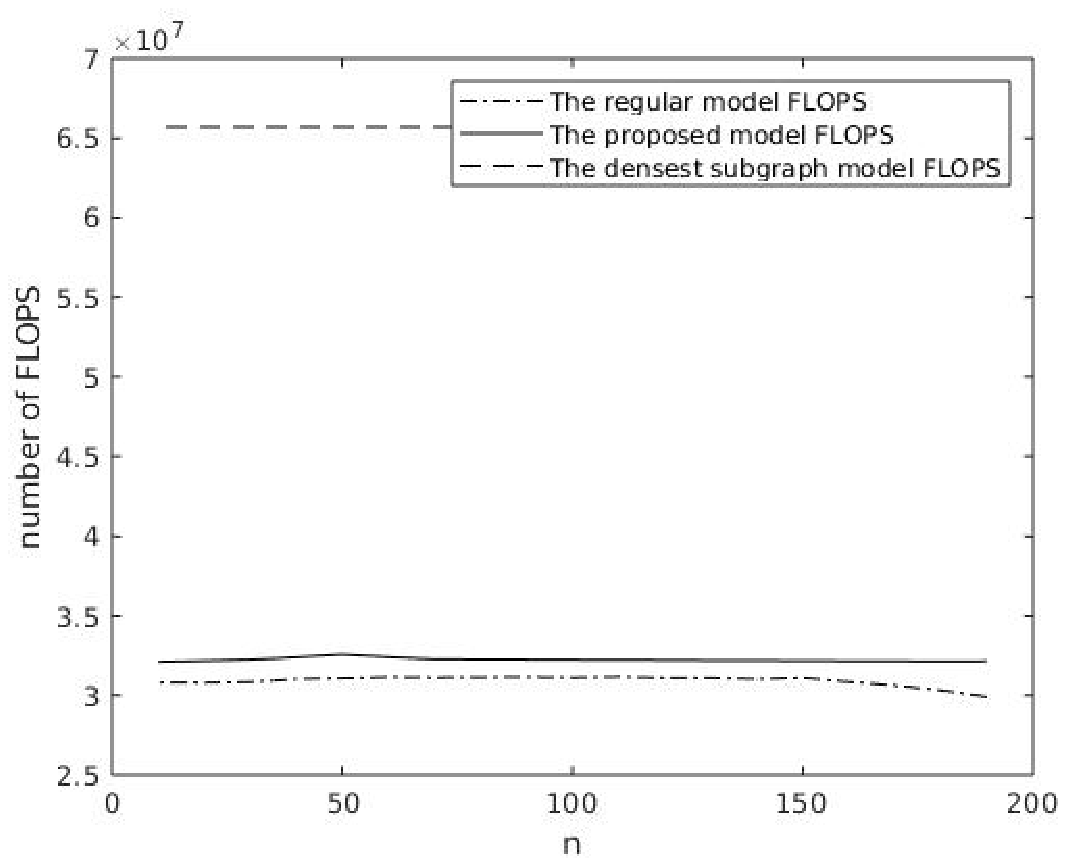}
		\label{fig:sub1}
	}
	\subfloat[Subfigure 2 list of figures text][$N=500$]
	{
		\includegraphics[width=0.5\textwidth]{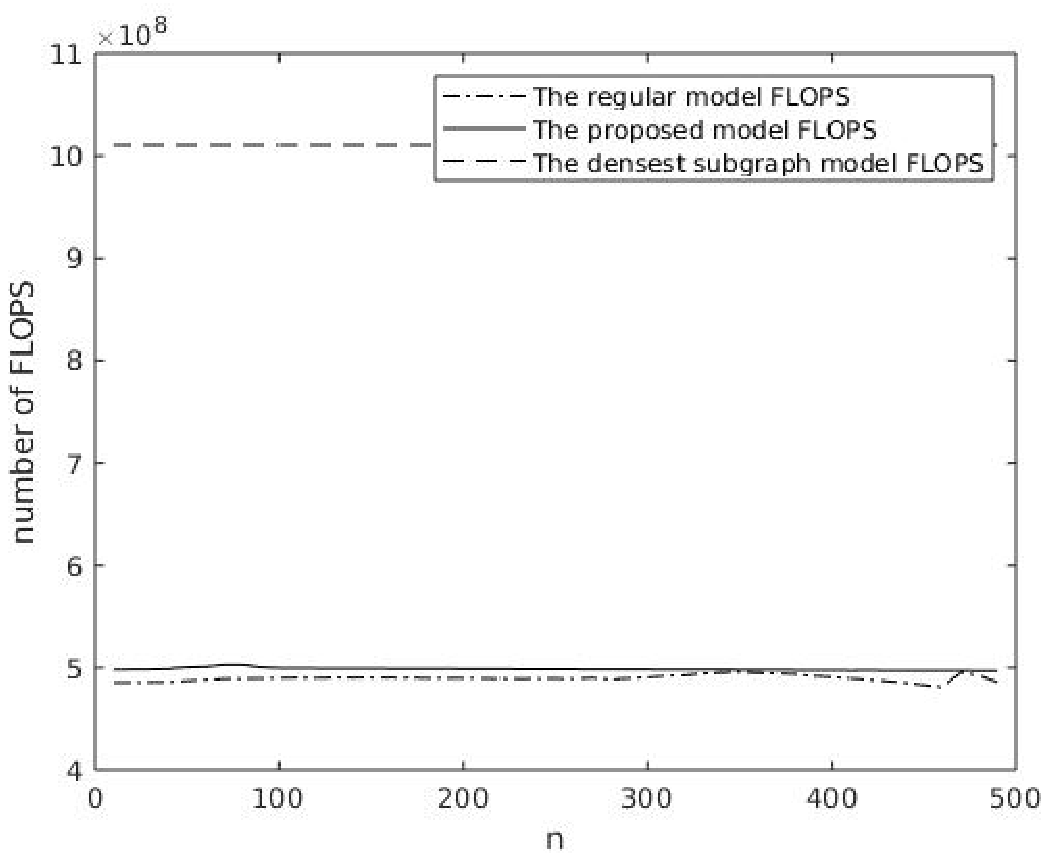}
		\label{fig:sub2}
	}
	\caption{Average number of FLOPS per iteration.}
	\label{fig:f}
\end{figure}

\vspace*{-5pt}
To clarify the scaling of our proposed approach, we sketch the number of FLOPS (Floating Point Operations per Second) needed per iteration. Figure \ref{fig:f} shows the average number of FLOPS needed per iteration for $N=200$ and $N=500$.

Figure \ref{fig:f} shows that the number of FLOPS needed per iteration is $\mathcal{O}(N^3)$. It also shows that our algorithm needed less number of FLOPS than the densest subgraph algorithm. 

Finally, we also compare our proposed model for the planted clique problem with two further known mathematical models. These are the nuclear norm minimization model (NNM) and the model based on semi-definite programming (SDP) \cite{ames2011convex}. NNM was solved using PPAPack, a software package in Matlab. NNM failed to obtain optimal solutions with the desired error tolerance better than $10^{-2}$ where the errors were found using \eqref{E1}  \cite{ames2011convex}. On the other hand, SDP solver failed to provide optimal planted clique of any size when $N\geq 100$. In addition, Ames \cite{ames2011convex} reported that final solutions of NNM had to be obtained by rounding the entries of solution matrix provided by the software, PPAPack, used. This is not the case for the optimal solutions obtained by our proposed algorithm, as we have claimed earlier in our paper, see Section \ref{Sec: Section 2}. In addition, the results presented in \cite{ames2011convex} show that all solutions were obtained with an error tolerance of $10^{-2}$ which is much inferior to our error tolerance of $10^{-8}$. Our proposed approach has probability 1 for all tested $(N,n)$ pairs, with $n\in [30,190]$, $n\in[50,490]$, and $n\in [60,990]$ for $N=200,~500$, and 1000, respectively. On the other hand NNM does not achieve probability 1 for all $(N,n)$ pairs with $n\in [110,140]$, $n\in[200,250]$, and $n\in [200,400]$ for $N=200,~500$, and 1000, respectively. Clearly, our algorithm solves harder problems than NNM and SDP.

\subsection{Maximum clique in random graphs}
We have also performed experiments on random graphs where all the edges are assigned with probability $p$, with no clique being planted. These results are presented in Table \ref{tab:Table11} where $n$ is the size of the maximum clique obtained by our algorithm. We have used the same stopping condition \eqref{SC1} to stop the algorithm, but calculated the errors using the formula 
\begin{align}\label{errorn}
Error = \Bigg   \vert \sqrt{\sum_{j=1}^{N} \sum_{i=1}^{N} L_{ij}} -\big \Vert L \big \Vert \Bigg  \vert,
\end{align}
since we have shown earlier in Section \ref{Sec:section 5n} that $\Vert L \Vert=n$. Our suggested formula \eqref{errorn} measures the recovery of maximum clique based on nodes in the clique, $\sqrt{\sum_{j=1}^{N} \sum_{i=1}^{N} L_{ij}}$.

The solution $(L,S)$ is initialized with $(0,0)$. A total of 32 runs have been performed for each $N$, and the results presented in each row of Table \ref{tab:Table11} are obtained for a single run. We have used $p=0.8$ and $p=0.87$ for $N=200$ and $N=500$ respectively. We have used higher values for $p$ to ascertain that a reasonable size of maximum clique is formed in each random graph. Results in Table \ref{tab:Table11} show the perfect recovery of maximum cliques based on nodes in the clique except for a small number of cases having some errors. These errors occur because the generated random graph has a bi-clique of size greater than the size of the clique. 
\vspace*{-5pt}
\begin{table}[ht]
	\begin{minipage}[b]{0.45\linewidth}\centering
		\setlength{\tabcolsep}{1pt}
		\begin{tabular}{|c|c|c|c|c|c|}
			\hline\noalign{\smallskip}
			\multicolumn{6}{|c|}{$N=200$}\\
			\hline 
			$n$ & $Error$ & Runtime&$n$ & $Error$ & Runtime\\
			\hline
			59 & 0 & 2.02 &46 & 0 & 1.97\\ 
			49 & 0 & 1.63 & 53 & 0 & 1.85\\ 
			74 & 0 & 1.99&57& 0 & 2.09\\ 
			59 & 0 & 1.88&43 & 0 & 1.75\\ 
			48 & 0 & 1.86&48 & 0 & 1.96\\ 
			55 & 0 & 1.86&51 & 0 & 1.89\\ 
			48 & 0 & 1.88&63 & 0 & 2.12\\ 
			54 & 0 & 1.82&58 & 0 & 1.88\\ 
			62 & 0 & 2.00&61 & 0 & 2.04\\ 
			57 & 0 & 2.10&50 & 0 & 1.84\\ 
			42 & 0 & 1.88&53 & 0 & 1.82\\ 
			47 & 0 & 2.35&44 & 0 & 1.92\\ 
			62.49 & 0.49 & 1.94&58& 0& 1.88\\ 
			62 & 0 & 1.76&42 & 0 & 1.83\\ 
			51 & 0 & 1.95&68 & 0 & 2.10\\ 
			45 & 0 & 1.83&36 & 0 & 1.87\\ 
			\noalign{\smallskip}\hline
		\end{tabular}
	\end{minipage}
	\hspace{0.5cm}
	\begin{minipage}[b]{0.45\linewidth}
		\setlength{\tabcolsep}{1pt}
		\begin{tabular}{|c|c|c|c|c|c|}
			\hline \noalign{\smallskip}
			\multicolumn{6}{|c|}{$N=500$}\\ 
			\hline 
			$n$ & $Error$ & Runtime&$n$ & $Error$ & Runtime\\
			\hline
			483 & 0 & 59.46&486.49& 0.49 & 54.77\\ 
			485.49 & 0.49 & 67.04&475 & 0 & 67.05\\ 
			482 & 0 & 65.00&494 & 0 & 54.23\\ 
			485 & 0 & 67.42&479 & 0 & 62.02\\ 
			481 & 0 & 67.35&483.49 & 0.49 & 62.03\\ 
			487 & 0 & 64.71&491 & 0 & 56.26\\ 
			485 & 0 & 57.99&481 & 0 & 57.16\\ 
			483.49 & 0.49 & 61.92&478 & 0 & 57.09\\ 
			491 & 0& 55.43&489.49& 0.49 & 54.71\\ 
			489 & 0 & 50.26&470 & 0 & 61.56\\ 
			490 & 0 & 63.37&491 & 0 & 55.39\\ 
			494 & 0 & 54.36&485 & 0 & 54.04\\ 
			495 & 0 & 56.35&488 & 0 & 53.82\\ 
			487 & 0 & 52.95&488 & 0 & 59.48\\ 
			494 & 0 & 54.23&488 & 0 & 63.48\\
			481 & 0 & 63.20&488 & 0 & 64.87\\ 
			\noalign{\smallskip}\hline
		\end{tabular}
	\end{minipage}
	\caption{Maximum clique in random graphs}
	\label{tab:Table11}
\end{table}
\vspace*{-5pt}
\subsection{Cliques in real-world graphs}
Our experiments include a few real-world graphs from the 10th DIMACs Implementation Challenge, which focus on clustering and partitioning graphs. The results of the real graphs are provided in Table \ref{tab:MCP-RWG}. We first consider the graph JAZZ, which is a representation of a collaboration network between Jazz musicians \cite{nr}. The nodes represent Jazz musicians, whereas the edges indicate that two musicians have collaborated in a band. The JAZZ graph consists of 198 vertices and 2742 edges. In a earlier study \cite{tsourakakis2013denser}, a clique of 30 vertices was found in this network. With the value $\rho = 0.25$, we employ \textbf{Algorithm} \ref{a4} in the adjacency matrix of this graph. After 37 iterations, our algorithm reaches the maximum clique of size 30 within 0.2152 second. We have implemented DSA with $\tau=0.35$, and stopped DSA with tolerance $10^{-4}$. It solves JAZZ in 0.5811 seconds with 94 iterations. We have also applied DSA to all 18 DIMACS benchmark problem and it failed in all problems.

%We now compare our algorithm with DSA using the problems in Table \ref{tab:MCP-RWG}. We first implement DSA on JAZZ using $n = 30$, $\tau= 0.35$ and stopping tolerance $10^{-4}$. The algorithm solves JAZZ in 0.5811 seconds with 94 iterations. DSA, however, fails to solve the other problems. 
We now compare our algorithm with the algorithm presented in \cite{belachew2017solving} using 18 DIMACS benchmark problems. We implement \textbf{Algorithm} \ref{a4}  with value $\rho=0.4$. Comparisons  are summarized in Table~\ref{tab:MCP-RWG}, where the symbol \lq -\rq ~ denotes non-availability of data. Results for the other algorithm under column 4, Table 2,  were taken from \cite{belachew2017solving}.

The number of iterations needed by \textbf{Algorithm} \ref{a4} is given in the last column.
Here $(N,\omega(G))$ represents the number of vertices and the clique number, respectively, while $n ( \textit{respectively}, n)$\cite{belachew2017solving} denotes the size of the clique obtained by our algorithm (respectively, by the algorithm in \cite{belachew2017solving}).
%We also consider some DIMACS benchmark data sets \cite{nr} and we compare our results with the results obtained in \cite{belachew2017solving}. We employ \textbf{Algorithm} \ref{a4}  with value $\rho=0.4$. Here $(N, \omega(G))$ represents the number of vertices and the clique number, respectively, while $n(\textit{respectively}, n\cite{belachew2017solving})$ denotes the size of the clique obtained by our algorithm (respectively, by the algorithm in \cite{belachew2017solving}). Results are summarized in Table \ref{tab:MCP-RWG}.
 
 \vspace*{-5pt}
 \begin{table}[H]
 	\setlength{\tabcolsep}{1pt}
 	\begin{tabular}{c c c c c c c c}
 		\hline
 		Graph & $(N,\omega(G))$\qquad \qquad &Number of edges \qquad \qquad & $n(n )$\cite{belachew2017solving} \qquad \qquad& Number of iterations\\
 		\hline \hline
 		%JAZZ& (198,30) & 2742 & 30(-) &37   \\
 		
 		BROCK200-1&  (200,21)&  14834 & 24 (19) &277  \\
 		
 		BROCK200-4&  (200,17)& 13089 & 34 (10)  & 124 \\
 		
 		BROCK400-2&  (400,29)& 59786 & 46 (24)  & 284 \\
 		
 		BROCK400-4&  (400,33)& 59765 & 37 (24)  & 324 \\
 		
 		C125.9 & (125,34) & 6963 & 34 (-) & 769\\
 		
 		C250.9 & (250, 44) & 27984 & 44(-) & 1012 \\
 		
 		C500.9 & (500,$\geq 57$)& 112332 & 216(50) & 1442 \\
 		C-fat500-10 & (500,-) & 46627 & 306(-) & 3 \\

 		GEN200-P0.9-44 & (200,44) & 17910 & 44 (-) & 1000 \\ 
 		
 		GEN200-P0.9-55 & (200,55) & 17910 & 55 (-) & 989 \\
 		
 		GEN400-P0.9-55 & (400,55) & 71820 & 134 (-) & 1442 \\ 
 		
 		GEN400-P0.9-65 & (400,65) & 71820 & 135 (-) & 1360 \\  
 		
 		GEN400-P0.9-75 & (400,75) & 71820 & 57 (-) & 1430 \\
 		
 		P-HAT300-2& (300,25) & 21928 & 40 (-) & 256  \\
 		P-HAT300-3& (300,36) & 33390 & 219 (-) & 225  \\
 		P-HAT500-2& (500,-) & 62946 & 159 (-) & 201  \\
 		P-HAT700-2& (700,44) & 121728 & 55 (-) & 478 \\
 		P-HAT700-3& (700,62) & 183010 & 209 (-) & 654  \\
 		\hline\hline
 	\end{tabular}
 	\\
 	\caption{Maximum cliques in real-world graphs}
 	\label{tab:MCP-RWG} 
 \end{table}

Comparison made in Table \ref{tab:MCP-RWG} shows that our  algorithm performs better than the algorithm in \cite{belachew2017solving} in the tested DIMACS benchmark data sets. Our algorithm recovers the confirmed clique sizes for 4  problems, while algorithm proposed in \cite{belachew2017solving}) failed to obtained confirmed clique for any problem. Moreover,
our proposed algorithm recovers cliques of large sizes as it decomposes the input
adjacency matrix of the input graph into a rank-one matrix and a sparse matrix. We have used our error formula in \eqref{errorn} to confirm the clique sizes recovered.  
%Comparison made in Table \ref{tab:MCP-RWG} shows that our proposed algprithm performs better than the algorithm in \cite{belachew2017solving} in the tested DIMACS benchmark data sets. Moreover, our proposed algorithm recovers the confirmed clique sizes for some problems, such as $C125.9$ and $C250.9$ and others.  

\section{Conclusions} \label{Sec: Section 6}
We have suggested a mathematical model for the clique problem that differs from the known matrix decomposition model in that it produces naturally integer solution required. This has been possible due to the dynamic nature of the weighted $\ell_{1}$-norm. We have established conditions that guarantee the recovery and the uniqueness of the solution, and we have derived a tight bound of the dual matrix that certifies the optimality conditions of our proposed model. Our approach produces much superior solution quality when compared to other known approaches. This has been possible due to the dynamic nature of our mathematical model. Our algorithm requires no input from the user other than the adjacency matrix of the input graph. In addition, the algorithm can be implemented at easy without needing any external solvers. Although the algorithm has been proposed for the planted clique problem, it has been tested on the maximum clique problem using random graphs with almost error-free results. We have also suggested a new expression for error calculations. Moreover, we have applied our algorithm to some real-world graphs and DIMACS data sets, and cliques are recovered successfully using our matrix decomposition model. 

\begin{acknowledgements}
This work is supported by the Organization for Women in Science from Developing World (OWSD) and Swedish International Development Cooperation Agency (Sida).
	The second author would like to thank Professor Stephen Vavasis of University of Waterloo for introducing him to the research topic.
\end{acknowledgements}

%\end{thebibliography}

\bibliographystyle{spmpsci} %%ieeetr
%\setcitestyle{square}
\bibliography{references}

\end{document}